\documentclass[12pt,draftcls, onecolumn]{IEEEtran}

\usepackage{cite}

\ifCLASSINFOpdf
  \usepackage[pdftex]{graphicx}
  \graphicspath{{./}}
  \DeclareGraphicsExtensions{.jpeg,.png}
\else
\fi
\usepackage{amsmath}
\usepackage{xcolor}
\usepackage{amsfonts}
\usepackage{amssymb}

\newtheorem{theorem}{Theorem}
\newtheorem{definition}{Definition}
\newtheorem{remark}{Remark}
\newtheorem{corollary}{Corollary}
\newtheorem{lemma}{Lemma}

\begin{document}
%
\title{Distribution of the Scaled Condition Number of  Single-spiked Complex Wishart Matrices}
%
%
%

\author{Pasan~Dissanayake,
  Prathapasinghe Dharmawansa,~\IEEEmembership{Member,~IEEE,}
  and~Yang~Chen
  \thanks{P. Dissanayake and P. Dharmawansa are with the Department
		of Electronic and Telecommunication Engineering, University of Moratuwa, Moratuwa 10400, Sri Lanka (e-mail: pasandissanayake@gmail.com, prathapa@uom.lk ).}%
  \thanks{Y. Chen is with the Department of Mathematics, Faculty of Science and Technology, University of Macau, Macau, P. R. China (e-mail:
 yayangchen@umac.mo).}
}
\maketitle

\begin{abstract}
Let $\mathbf{X}\in\mathbb{C}^{n\times m}$ ($m\geq n$) be a random matrix with independent columns each distributed as complex multivariate Gaussian with zero mean and {\it single-spiked} covariance matrix $\mathbf{I}_n+ \eta \mathbf{u}\mathbf{u}^*$, where $\mathbf{I}_n$ is the $n\times n$ identity matrix, {\color{blue}$\mathbf{u}\in\mathbb{C}^{n\times 1}$} is an arbitrary vector with unit Euclidean norm, $\eta\geq 0$ is a non-random parameter, and $(\cdot)^*$ represents the conjugate-transpose. This paper investigates the distribution of the random quantity $\kappa_{\text{SC}}^2(\mathbf{X})=\sum_{k=1}^n \lambda_k/\lambda_1$, where {\color{blue}$0\le \lambda_1\le \lambda_2\le \ldots\leq \lambda_n<\infty$} are the ordered eigenvalues of $\mathbf{X}\mathbf{X}^*$ (i.e., single-spiked Wishart matrix). This random quantity is intimately related to the so called {\it scaled condition number} or the Demmel condition number (i.e., $\kappa_{\text{SC}}(\mathbf{X})$) and the minimum eigenvalue of the fixed trace Wishart-Laguerre ensemble (i.e., $\kappa_{\text{SC}}^{-2}(\mathbf{X})$). In particular, we use an orthogonal polynomial approach to derive an exact expression for the probability density function of $\kappa_{\text{SC}}^2(\mathbf{X})$ which is amenable to asymptotic analysis as matrix dimensions grow large. Our asymptotic results reveal that, as $m,n\to\infty$ such that $m-n$ is fixed and when $\eta$ scales on the order of $1/n$, $\kappa_{\text{SC}}^2(\mathbf{X})$ scales on the order of $n^3$. In this respect we establish simple closed-form expressions for the limiting distributions.  {\color{blue}It turns out that, as $m,n\to\infty$ such that $n/m\to c\in(0,1)$, properly centered $\kappa_{\text{SC}}^{2}(\mathbf{X})$ fluctuates on the scale $m^{\frac{1}{3}}$}.

\end{abstract}

\begin{IEEEkeywords}
  Condition number, cumulative distribution function (c.d.f.), eigenvalues, hypergeometric function of two matrix arguments, moment generating function (m.g.f.), orthogonal polynomials, probability density function (p.d.f.), single-spiked covariance, Wishart matrix. 
\end{IEEEkeywords}

%
\IEEEpeerreviewmaketitle

\section{Introduction}
\IEEEPARstart{C}{ondition} numbers measure the worst-case sensitivity of problems with respect to small perturbations of the input. 
The seminal studies by Turing \cite{ref:Turing} and John von Neumann and Herman Goldstine \cite{ref:neumannGoldstine} had introduced the condition number as a powerful tool to understand the loss of precision in the solution of linear systems in the presence of finite-precision arithmetic. 
This particular tool is fundamentally important in numerical linear algebra and related areas.
However, as conjectured in \cite{ref:renegar}, computing the condition number corresponding to a certain input, for a given problem, is no easier than solving the problem itself for this particular input. To circumvent this difficulty, the concept of probabilistic analysis of condition numbers has been introduced; see e.g., \cite{ref:cucker} for a partial list of references in this respect. The key concept here is to endow the input set with a certain probability measure and thereby statistically characterize the condition number. Since this characterization rely on the underlying probability measure, the current literature identifies two different approaches, namely, {\it average analysis} \cite{ref:demmelProbability} and {\it smoothed analysis} \cite{ref:spielman, ref:burgisser,ref:burgisser2,ref:burgisser3} depending on the nature of the input distribution. In particular, former approach assigns real/complex standard normal density to the input set, whereas the latter approach assumes that the input set is endowed with real/complex non-zero mean normal measure.

Among various condition numbers, perhaps the best known condition number with respect to matrix inversion problem assumes \cite{ref:neumannGoldstine} $\kappa(\mathbf{A})=||\mathbf{A}||_2 ||\mathbf{A}^{-1}||_2$, where $\mathbf{A}\in\mathbb{C}^{n\times n}$ and $||\cdot||_2$ denotes the $2-$norm. The statistical characteristics of $\kappa(\mathbf{A})$ have been well studied in the literature \cite{ref:demmel, ref:kostlan, ref:smale, ref:edelman}. Another fundamental form introduced by Demmel in his seminal work on the probabilistic analysis of the degree of difficulty associated with numerical analysis problems \cite{ref:demmelProbability} is defined as $\kappa_D(\mathbf{A})=||\mathbf{A}||_F ||\mathbf{A}^{-1}||_2$, where $||\cdot||_F$ denotes the Frobenius norm. This definition naturally extends to rectangular matrices \cite{ref:prathapSIAM,ref:burgisser3} giving
\begin{align}
    \kappa_{\text{SC}}(\mathbf{X})=||\mathbf{X}||_F ||\mathbf{X}^\dagger||_2
\end{align}
where $\mathbf{X}\in \mathbb{C}^{n\times m}$ and $(\cdot)^\dagger$ denotes the Moore-Penrose or pseudo inverse operation \cite{ref:golub}. We refer to  $\kappa_{\text{SC}}(\mathbf{X})$ as the scaled condition number (SCN)\footnote{This belongs to a certain class known as the {\it conic} condition numbers \cite{ref:burgisser3}.}; a term coined by Edelman in \cite{ref:edelmanSCN}. A convenient alternative representation of $\kappa_{\text{SC}}(\mathbf{X})$ involving the spectral characteristics of $\mathbf{X}$ can be written as $
    \kappa_{\text{SC}}(\mathbf{X})=\sqrt{\sum_{k=1}^r \lambda_k/\lambda_1}
$,
where $\text{rank}(\mathbf{X})=r\leq \min (m,n)$ and $\lambda_1\leq \lambda_2\leq \ldots\leq \lambda_r$ are the non-zero eigenvalues of $\mathbf{X}\mathbf{X}^*$ (or $\mathbf{X}^*\mathbf{X}$) with $(\cdot)^*$ denoting the conjugate transpose operator. Since the statistical characteristics of $\kappa_{\text{SC}}(\mathbf{X})$  are of paramount importance in many scientific disciplines, it is common to assume $\mathbf{X}$ to be real/complex Gaussian distributed with $m\geq n$ \cite{ref:prathapSIAM,ref:edelmanSCN,ref:edelman,ref:wei,ref:zhong} which in turn gives
\begin{align}
    \kappa_{\text{SC}}(\mathbf{X})=\sqrt{\frac{\sum_{k=1}^n \lambda_k}{\lambda_1}}.
\end{align}

The statistical characteristics of the SCN and  its variants\footnote{The quantity $\kappa^{-2}_{\text{SC}}(\mathbf{X})$ is known as the minimum eigenvalue of the fixed trace Wishart-Laguerre ensemble in the statistical physics literature; see e.g \cite{ref:dean,ref:vivo,ref:mehta,ref:forresterLogGases,ref:kumar,ref:chen,ref:akemann} and references therein.} have been instrumental in understanding many physical phenomena across a heterogeneous fields of sciences \cite{ref:demmelProbability,ref:burgisser3,ref:edelmanSCN,ref:prathapSIAM,ref:prathapJMVA,ref:heath,ref:kobyakov,ref:liJianzhiClusterBased,ref:liJianzhi3DCluster,ref:caiJia,ref:krishnaiahPart1,ref:krishnaiahPart2,ref:krishnaiahEvaluation,ref:davis,ref:dean,ref:vivo,ref:mehta,ref:forresterLogGases,ref:kumar,ref:chen}. While numerical analysts and statistical physicists are interested in the behavior of $\kappa_{\text{SC}}(\mathbf{X})$ for $\mathbf{X}$ having independent complex normal entries, the case corresponding to {\it correlated} complex normal entries are instrumental in wireless communications and statistics \cite{ref:heath,ref:kobyakov,ref:liJianzhiClusterBased,ref:liJianzhi3DCluster,ref:caiJia,ref:krishnaiahPart1,ref:krishnaiahPart2,ref:krishnaiahEvaluation}. In particular,  $\kappa_{\text{SC}}(\mathbf{X})$ has been used as a performance metric in certain wireless signal processing applications involving multiple-input multiple-output (MIMO) systems \cite{ref:heath,ref:liJianzhiClusterBased,ref:liJianzhi3DCluster} in which $\mathbf{X}$ corresponds to the rich scattering matrix channel between the transmit and receive antenna arrays \cite{ref:tse}. Moreover, in the presence of antenna correlation, $\mathbf{X}$ is commonly modelled as a {\it correlated} Gaussian random matrix; e.g, see \cite{ref:chiani} for a partial list of references. Therefore, these facts further highlight the utility of
$\kappa_{\text{SC}}(\mathbf{X})$ for random and correlated $\mathbf{X}$.
The exact statistical characteristics of $\kappa_{\text{SC}}(\mathbf{X})$ for $\mathbf{X}$ having independent real/complex entries are well documented in the literature \cite{ref:wei,ref:zhong,ref:prathapSIAM,ref:edelmanSCN,ref:vivo,ref:mehta,ref:forresterLogGases,ref:kumar,ref:chen,ref:akemann}. 

{\color{blue} Among various covariance structures, Johnstone's  spiked model \cite{ref:johnstone} has been widely used in the literature to analyze the effects of having a few dominant trends or correlations in the covariance matrix. To be precise, under this setting, the covaraince matrix $\boldsymbol{\Sigma}\in\mathbb{C}^{n\times n}$ of $\mathbf{X}$ is modeled as $\boldsymbol{\Sigma}=\mathbf{I}_n+\sum_{k=1}^r \theta_k\mathbf{u}_k\mathbf{u}_k^*$, where $\mathbf{I}_n$ is the $n\times n$ identity matrix,  $\mathbf{u}_k\in \mathbb{C}^{n\times 1}$, $k=1,2,\ldots,r (\leq n)$ are a set of orthonormal vectors, and $\theta_1\geq \theta_2\geq \ldots\geq \theta_r\geq 0$. 
Consequently, the $\mathbf{u}_k$s' are referred to as the {\it spikes} and this particular covariance structure is sometimes known as rank-$r$ perturbation of the identity matrix. This fact is further highlighted by the eigen-structure of $\boldsymbol{\Sigma}$
in which the the dominant $r$ eigenvalues  can be written as, $\theta_1+1\geq \theta_2+1\geq \ldots\geq \theta_r+1$, whereas the rest of the $n-r$ eigenvalues assume $1$.
These spikes arise in various practical scenarios in different scientific disciplines. For instance, they correspond to the first few dominant factors in factor models arising from financial economics \cite{ref:fan,ref:ontaski,ref:torun}, first few principal components \cite{ref:perry,ref:johnstone}, the number of clusters in gene expression data \cite{ref:ke}, and the number of signals in signal detection and estimation, see e.g., \cite{ref:couillet,ref:onatskiSignal,ref:hallin,ref:raoTIT,ref:najim,ref:hachem,ref:chamain,ref:prathapRoyRoot,ref:baikPhaseTrans, ref:montanari,ref:haddadi}, and references therein. In particular, \cite{ref:onatskiSignal} and \cite{ref:raoTIT} have focus on the rank-$r$ model, whereas the rank-$1$ (i.e., single-spiked) model, which is of our interest in this manuscript, has been employed by \cite{ref:najim,ref:hachem,ref:chamain,ref:prathapRoyRoot,ref:baikPhaseTrans, ref:montanari,ref:haddadi,ref:boaz,ref:hallin} in the signal detection problem. In a sharp contrast, Hanlen and Grant \cite{ref:hanlen} have used the rank$-r$ model to investigate the effect of correlation on the MIMO capacity. Be that as it may, the constant correlation model, which is one of the most important correlation models frequently used in a wide array of MIMO applications \cite{ref:ping,ref:shin,ref:win,ref:marco,ref:simon,ref:biglieri}, gives rise to a single spiked model for the scaled $\boldsymbol{\Sigma}$ matrix. To be specific, under this setting, $\boldsymbol{\Sigma}$ consists of $1$'s in the main diagonal and $\sigma\in[0,1)$'s in all off-diagonal entries. As such, $\frac{1}{1-\sigma}\boldsymbol{\Sigma}$ admits the desired single spiked structure given by $\frac{1}{1-\sigma}\boldsymbol{\Sigma}=\mathbf{I}_n+\frac{n\sigma}{1-\sigma} \mathbf{11}^*$, where $\mathbf{1}=\left(\frac{1}{\sqrt{n}}\;\frac{1}{\sqrt{n}}\;\ldots\; \frac{1}{\sqrt{n}}\right)^*$. For instance, this  has been exploited in \cite{ref:ping,ref:shin} to  derive certain performance measures related to MIMO systems. Moreover, the SCN has been  used as a performance metric in several wireless signal processing applications involving MIMO systems as delineated in \cite{ref:heath,ref:liJianzhiClusterBased,ref:liJianzhi3DCluster}. Therefore, these facts further highlight the utility of the single spiked covariance model in a wide class of applications. 
 
 The square of the SCN $\kappa^2_{\text{SC}}(\mathbf{X})$ in conjunction with correlated Gaussian $\mathbf{X}$ having a single-spiked covariance structure has been instrumental in the spectrum sensing problem of cognitive radio (CR) networks\footnote{The key concept behind CR is to opportunistically utilize the underutilized spectrum in view of improving the spectral efficiency \cite{ref:mitola,ref:goldsmith,ref:larson}.} \cite{ref:liang,ref:larson,ref:ayse}. In particular, one of the blind detection techniques proposed in the seminal paper \cite{ref:liang} uses $\kappa^2_{\text{SC}}(\mathbf{X})$ as the test statistic to detect the presence of primary user (i.e., to detect whether the particular slice of spectrum is occupied by the intended user or not). To be specific, let us consider a scenario where the secondary user is equipped with $n>1$ antennas (or sensors). Then the spectrum sensing problem can be formulated as the following binary hypothesis testing problem \cite{ref:larson,ref:najim}:
 \begin{align}
     \mathcal{H}_1: \mathbf{x}(k)&=\mathbf{h}s(k)+\mathbf{w}(k),\;\; k=1,2,\ldots,m\nonumber\\
     \mathcal{H}_0: \mathbf{x}(k)&=\mathbf{w}(k),\;\; k=1,2,\ldots,m
 \end{align}
 where $\mathcal{H}_1$ and $\mathcal{H}_0$ are the ``primary signal present" and ``primary signal absent" hypotheses, $\mathbf{x}(k)\in\mathbb{C}^{n\times 1}$ is the observed signal, $\mathbf{h}\in \mathbb{C}^{n\times 1}$ is the channel between the source and the detector, $s(k)\sim\mathcal{CN}(0,\gamma)$ denotes the transmitted signal with $\gamma>0$, $\mathbf{w}(k)\sim\mathcal{CN}_{n}(\mathbf{0},\sigma^2 \mathbf{I}_n)$ denotes the noise process, and $m\geq n$ denotes the number of {\it independent} observations (samples). Consequently, the population covariance matrix can be written as
 \begin{align}
 \label{hypotheses}
     \mathbf{R}=\mathbb{E}\left\{\mathbf{x}(k)\mathbf{x}(k)^*\right\}=\left\{\begin{array}{ll}
     \gamma \mathbf{hh}^* +\sigma^2 \mathbf{I}_n & \text{under $\mathcal{H}_1$}\\
     \sigma^2\mathbf{I}_n & \text{under $\mathcal{H}_0$}
     \end{array}
     \right.
 \end{align}
 where $\mathbb{E}\left\{\cdot\right\}$ denotes the mathematical expectation operator. Now, in the absence of the knowledge of the parameters $\mathbf{h},\gamma$, and $\sigma^2$ at the receiver, one of the blind-eigenvalue based test statistics proposed in \cite{ref:liang} (also in \cite{ref:larson} and \cite{ref:ayse}) assumes (see e.g., \cite{ref:larson} for a comprehensive discussion on this and related other detectors)
 \begin{align}
     T(\hat{\boldsymbol{\lambda}})= \frac{ \sum_{k=1}^n \hat{\lambda}_k}{\hat{\lambda}_1}
 \end{align}
 where $\hat{\lambda}_1\le \hat{\lambda_2}\le \ldots\le \hat{\lambda}_n$ are the ordered eigenvalues of the sample covariance matrix
 \begin{align}
     \hat{\mathbf{R}}=\frac{1}{m}\sum_{k=1}^m \mathbf{x}(k)\mathbf{x}(k)^*
 \end{align}
 and $\hat{\boldsymbol{\lambda}}=(\hat{\lambda}_1\; \hat{\lambda_2}\;\ldots\; \hat{\lambda}_n)$. Since $T(\hat{\boldsymbol{\lambda}})$ is not affected if $\hat{\mathbf{R}}$ is scaled by a given constant, we find it convenient to rewrite the above test statistic as
 \begin{align}
 \label{Teststat}
     T(\boldsymbol{\lambda})=\kappa_{\text{SC}}^2(\mathbf{X})=\frac{ \sum_{k=1}^n {\lambda}_k}{{\lambda}_1}
 \end{align}
 where $\lambda_1\le \lambda_2\le \ldots\le \lambda_n$ are the ordered eigenvalues of $\mathbf{XX}^*$ with $\mathbf{X}=\left(\mathbf{x}(1) \; \ldots \; \mathbf{x}(m)\right)\in\mathbb{C}^{n\times m}$. Therefore, in view of (\ref{hypotheses}) and noting that the columns of $\mathbf{X}$ are independent, we obtain 
 \begin{align}
 \label{CRsetting}
     \mathbf{X}\sim\left\{\begin{array}{cc}
     \mathcal{CN}_{n,m}\left(\mathbf{0},\boldsymbol{\Sigma}\otimes \mathbf{I}_m\right) & \text{under $\mathcal{H}_1$}\\
     \mathcal{CN}_{n,m}\left(\mathbf{0},\mathbf{I}_n\otimes \mathbf{I}_m\right) & \text{under $\mathcal{H}_0$}
     \end{array}
     \right.
 \end{align}
 where $\boldsymbol{\Sigma}=\mathbf{I}_n+\frac{\gamma||\mathbf{h}||^2}{\sigma^2} \mathbf{u}_{h} \mathbf{u}_h^*$ denotes the single-spiked covariance model with $\mathbf{u}_h=\mathbf{h}/||\mathbf{h}||$, $\otimes$ stands for the Kronecker product of two matrices, and $||\cdot||$ denotes the Euclidean norm. It is  noteworthy that the term $\displaystyle \frac{\gamma ||\mathbf{h}||^2}{\sigma^2}$ is also known as the signal-to-noise ratio (SNR) \cite{ref:najim}. Clearly, to study the power of the above test for a given false alarm rate, one needs to statistically characterize the density of $T(\boldsymbol{\lambda})|\mathcal{H}_1$, thereby the density of $\kappa_{\text{SC}}^2(\mathbf{X})|\mathcal{H}_1$. Such a characterization has been confined to $n=2$ scenario only\cite{ref:matth}. This further highlights the utility of the  finite dimensional statistical characterization of the distribution of $\kappa_{\text{SC}}^2(\mathbf{X})$ for $\mathbf{X}$ with a single-spiked covariance structure.

Having motivated with the above facts, as the main contribution of this paper, we address the generic problem of determining the probability density function (p.d.f.) of $\kappa_{\text{SC}}^2(\mathbf{X})$ for  $\mathbf{X}\sim \mathcal{CN}_{n,m}\left(\mathbf{0},\boldsymbol{\Sigma}\otimes\mathbf{I}_m\right)$ with $\boldsymbol{\Sigma}=\mathbf{I}_n+\eta \mathbf{uu}^*$, where $\eta\geq 0$ and $\mathbf{u}\in\mathbb{C}^{n\times 1}$ is a unit vector. }
In particular, here we adopt a moment generating function (m.g.f.) based approach to derive a novel expression for the p.d.f. of $\kappa^2_{\text{SC}}(\mathbf{X})$.
The single-spiked covariance structure in turn enables us to leverage the powerful orthogonal polynomial techniques developed in \cite{ref:mehta} to arrive at our final p.d.f. expression. This novel expression developed for the p.d.f. of $\kappa^2_{\text{SC}}(\mathbf{X})$ contains the determinant of a square matrix whose dimensions depend on the relative difference between $m$ and $n$ (i.e., $m-n$). {\color{blue} For instance, in the CR setting discussed above, this refers to the relative difference between the number of observed samples and the number of secondary user antennas (or sensors). Moreover, this new p.d.f. expression further facilitates the evaluation of the power (i.e., probability of detection) of the test $T$. To further highlight this fact, we generate receiver operating characteristic (ROC) curves of the test $T$ for various $n,m$ configurations. Although obtaining an explicit functional relationship between the detection and false alarm probabilities seems an arduous task, for the important configuration of $m=n$, an explicit relation between those quantities has been derived.  It turns out that increasing either the SNR or the relative difference between $m$ and $n$ leads to an improved detection probability.}

{\color{blue}Whereas the above statistical characterization is valid for arbitrary $m$ and $n$, it is of paramount importance to investigate the behavior of $\kappa^2_{\text{SC}}(\mathbf{X})$ (or $T(\boldsymbol{\lambda})$) in various asymptotic domains. This amounts to establishing stochastic convergence result for the random variable $\kappa^2_{\text{SC}}(\mathbf{X})$. In this respect, it is critical to choose the correct asymptotic domain. For instance, related to multi-antenna communications, it is common to consider the domain in which the number of transmit and receive antennas diverge to infinity such that their ratio is fixed, see e.g.,  \cite{ref:couillet,ref:tulino} and references therein. Although this ratio can be any non-negative real number, it is taken to be {\it unity} in the analysis of asymptotic capacity scaling limits in multi-antenna systems \cite{ref:tseChuah,ref:liu,ref:akbar}. Against this backdrop, noting that the algebraic complexity of the new p.d.f. depends on $m-n$, we choose the asymptotic domain in which $m,n\to\infty$ such that $m-n$ is fixed. To be specific, we assume that\footnote{Given the functions $f(x)$ and $g(x)$, we have, for $x\to\infty$, $f(x) = O(g(x))$ if $\limsup_{x\to\infty} |f(x)/g(x)| < \infty$. } $m-n=O(1)$. This particular scaling has been instrumental in applications involving signal detection \cite{ref:chamain}, capacity scaling in multi-antenna systems  \cite{ref:tseChuah,ref:liu,ref:akbar}, multivariate analysis \cite{ref:prathapJMVA}, and  theoretical physics \cite{ref:vivo,ref:chen,ref:Forrester,ref:kumar,ref:tim}. In the setting of CR spectrum sensing, this is tantamount to the assumption that the number of samples (i.e., $m$) and the secondary user antennas (i.e., $n$) are of the same order but diverging. However, in practice, the number of co-located antennas at the secondary user cannot be increased due to space limitations. This drawback can be alleviated if we consider either a geographically distributed array of antennas at the secondary user or geographically distributed multiple secondary users each with a single antenna \cite{ref:ratna}\footnote{This is also known as cooperative spectrum sensing \cite{ref:larson,ref:ratna}.}. Therefore, the above facts further accentuate the utility of the asymptotic characterization of $\kappa^2_{\text{SC}}(\mathbf{X})$ in the domain $m-n=O(1)$. Capitalizing on this, we prove that the scaled random variable $\kappa^2_{\text{SC}}(\mathbf{X})/n^3$ converges in distribution to a random variable whose cumulative distribution function (c.d.f.) as well as p.d.f. contain the Bessel kernels. These limiting distributions are compact and having less algebraic complexity than their finite dimensional counterparts. Our numerical results reveal that, although derived for asymptotically large $m$ and $n$, these distributions compare favorably with finite values of those parameters as well. The limiting c.d.f. expression has been derived based on a new compact c.d.f. that we have obtained for the minimum eigenvalue of the single-spiked Wishart matrix. Although various c.d.f. and p.d.f. expressions have been derived in the literature for the minimum eigenvalue of this particular ensemble \cite{ref:tim,ref:petermin,ref:rathnasiam,ref:zanella,ref:chianiIT,ref:chianiRMT}, our newly derived expression is more compact and algebraically less complicated than those expressions. Apart from this, we also derive a new p.d.f. for the minimum eigenvalue as well. It turns out  that this expression is more simpler than that one can derive by taking the derivative of the corresponding c.d.f.

Another pertinent asymptotic domain is characterized by $m,n\to\infty$ such that $n/m\to c\in(0,1)$, see e.g., \cite{ref:tulino, ref:najim, ref:baikPhaseTrans, ref:johnstone,ref:onatskiSignal} and references therein. In this regime, stochastic convergence result for $\kappa^{2}_{\text{SC}}(\mathbf{X})$ has been established. In particular, we show that, for $\eta=O(1)$, properly centered and scaled $\kappa^{2}_{\text{SC}}(\mathbf{X})$ converges in law to the famous Tracy-Widom distribution \cite{ref:tracy} corresponding to $\beta=2$ (i.e., complex case)\footnote{Here $\beta>0$ is a non-random parameter which assumes $1,2,4$ for real ($\mathbb{R}$) symmetric, complex ($\mathbb{C}$) Hermitian, and quaternion ($\mathbb{H}$) self-dual case, respectively.}. This further reveals that suitably centered $\kappa^{2}_{\text{SC}}(\mathbf{X})$ fluctuates on the scale $m^{1/3}$. Since the above limiting distribution encompasses the case $\eta=0$ (i.e., when $\boldsymbol{\Sigma}=\mathbf{I}_n$ or equivalently no primary user signal is present) as well, we conclude that $\kappa^{2}_{\text{SC}}(\mathbf{X})$ (also $T(\boldsymbol{\lambda})$) does not have statistical power to detect a weak signal. Nevertheless,  this conclusion may not necessarily be true in the presence of a {\it strong} primary user signal. Moreover, we have shown that properly centered and scaled random variable $\kappa^{-2}_{\text{SC}}(\mathbf{X})$ also converges in law to the same Tracy-Widom distribution in this asymptotic regime. However, suitably centered $\kappa^{-2}_{\text{SC}}(\mathbf{X})$ fluctuates on the scale $m^{-5/3}$.}



The remainder of this paper is organized as follows. Section II provides some key preliminary results required in the subsequent sections. The new exact p.d.f. of $\kappa^2_{\text{SC}}(\mathbf{X})$ is derived in Section III. It also gives certain particularizations of the general p.d.f. expression. Apart from these results, we present the ROC curves corresponding to the test $T$ in the above outlined CR setting.  Moreover, a detailed asymptotic analysis of   $\kappa^2_{\text{SC}}(\mathbf{X})$ is provided in Section IV. Finally, conclusive remarks are made in Section V.

\section{Preliminaries}
To facilitate our main derivations, we will require the following preliminary results and definitions.
\begin{definition}
  Let $\mathbf{X}\in\mathbb{C}^{n\times m}$ ($m\geq
  n$) be distributed as $ \mathcal{CN}_{n,m}\left(\mathbf{0},\boldsymbol{\Sigma}\otimes\mathbf{I}_m\right)$, where $\boldsymbol{\Sigma}\in\mathbb{C}^{n\times n}$ is a Hermitian positive definite matrix. Then the matrix $\mathbf{W}=\mathbf{X} \mathbf{X}^*$ is said to follow a complex correlated Wishart distribution, i.e., $\mathbf{W}\sim\mathcal{W}_n(m,\boldsymbol{\Sigma})$.
\end{definition}

\begin{theorem}
  \label{thm1}
  The joint density of the ordered eigenvalues $0<\lambda_1\leq\cdots\leq\lambda_n$ of $\mathbf{W}$ is given by \cite{ref:james}
  \begin{align}
    \label{eq_joint_pdf_ordered_eigenvals}
    f\left(\lambda_1,\lambda_2,\ldots,\lambda_n\right)=\frac{K_{n,\alpha}}{\det^m(\boldsymbol{\Sigma})} \prod_{i=1}^n\lambda_i^\alpha \Delta_n^2(\boldsymbol{\lambda})
    {}_0\widetilde{F}_{0}\left(\boldsymbol{\Lambda},-\boldsymbol{\Sigma}^{-1}\right)
  \end{align}
  where 
  \begin{equation*}
    K_{n,\alpha}=\frac{1}{\prod_{i=1}^n(n+\alpha-i)!(n-i)!},
  \end{equation*}
  $\boldsymbol{\Lambda}=\rm{diag}(\boldsymbol{\lambda})$ with $\boldsymbol{\lambda}=\left(\lambda_1,\ldots,\lambda_n\right)$,  $\Delta_n(\boldsymbol{\lambda})=\prod_{1\leq i<k\leq n}\left(\lambda_k-\lambda_i\right)$ denotes the Vandermonde determinant, ${}_0\widetilde F_0(\cdot;\cdot)$ denotes the complex hypergeometric function of two matrix arguments, $\rm{diag}(\cdot)$ denotes the diagonal matrix, and $\text{det}(\boldsymbol{\Sigma})$ is the determinant of the square matrix $\boldsymbol{\Sigma}$. For $n\times n$ Hermitian  matrices $\mathbf{S}$ and $\mathbf{T}$, we have \cite{ref:james}
  \begin{align*}
    {}_0\widetilde F_0\left(\mathbf{S},\mathbf{T}\right)=\sum_{k=0}^\infty \frac{1}{k!} \sum_{\kappa} \frac{\mathcal{C}_\kappa(\mathbf{S})\mathcal{C}_\kappa(\mathbf{T})}{ \mathcal{C}_\kappa(\mathbf{I}_n)}
  \end{align*} 
  where $\mathcal{C}_\kappa(\cdot)$ is the complex zonal polynomial\footnote{The zonal polynomial $\mathcal{C}_\kappa(\mathbf{A})$ is a symmetric, homogeneous polynomial of degree $k$ in the eigenvalues of $\mathbf{A}$. However, the exact definition of the zonal polynomial is tacitly avoided, since it is not required in the subsequent analysis. More details of the zonal polynomials can be found in \cite{ref:james,ref:takemura}.}, $\kappa=(k_1,\ldots,k_n)$, with $k_i$'s being non-negative integers, is a partition of $k$ such that $k_1\geq\cdots\geq k_n\geq 0$ and $\sum_{i=1}^nk_i=k$. Moreover, ${}_0\widetilde F_0(\mathbf{S};\mathbf{T})$ admits the following unitary integral form \cite{ref:james}
  \begin{align}
  \label{uniint}
      {}_0\widetilde F_0(\mathbf{S};\mathbf{T})=\int_{\mathcal{U}_n} e^{\text{tr}\left(\mathbf{S}\mathbf{U}\mathbf{T}\mathbf{U}^*\right)} {\rm d}\mathbf{U}
  \end{align}
  where $\mathbf{U}\in \mathbb{C}^{n\times n}$ is a unitary matrix, ${\rm d}\mathbf{U}$ denotes the invariant measure (i.e., Haar measure) on the unitary group $\mathcal{U}_n$ normalized to make the total measure one (i.e., $\int_{\mathcal{U}_n} {\rm d}\mathbf{U}=1$), and  $\text{tr}(\cdot)$ is the trace of a square matrix.
\end{theorem}

\begin{remark}
Alternatively, following \cite{ref:khatri}, we have the following determinant representation of the complex hypergeometric function of two matrix arguments
\begin{align}
    {}_0\widetilde F_0\left(\mathbf{S},\mathbf{T}\right)=
    \prod_{k=1}^n (n-k)!\frac{\det\left[e^{s_it_j}\right]_{i,j=1,\ldots,n}}{\Delta_n(\mathbf{s})\Delta_n(\mathbf{t})}
\end{align}
where $\mathbf{s}=\{s_1,s_2,\ldots,s_n\}$, $\mathbf{t}=\{t_1,t_2,\ldots,t_n\}$ are the eigenvalues of $\mathbf{S}$ and $\mathbf{T}$ respectively, and the determinant of an $n\times n$ matrix with the $(i,j)$\textsuperscript{th} element given by $e^{s_it_j}$ is denoted by $\det\left[e^{s_it_j}\right]_{i,j=1,\ldots,n}$.
\end{remark}

Since we are interested in the so called {\it single-spiked} covariance model for $\boldsymbol{\Sigma}$, {\color{blue}following \cite{ref:wang,ref:peterBeta}, the joint eigenvalue density given in Theorem \ref{thm1} can be simplified to yield the expression given in the following corollary.}
\begin{corollary}
  Let $\boldsymbol{\Sigma}=\mathbf{I}_n+\eta \mathbf{u}\mathbf{u}^*$, where $\mathbf{u}\in\mathbb{C}^{n\times 1}$ with $||\mathbf{u}||=1$, and $\eta\geq 0$. Then the joint eigenvalue density of $\mathbf{W}\sim\mathcal{W}_n(m,\boldsymbol{\Sigma})$ (i.e., single-spiked Wishart-Laguerre ensemble\footnote{Alternatively, we can refer to it as the deformed Wishart-Laguerre ensemble.}) assumes
  \begin{align}
    \label{cor joint}
    f(\lambda_1,..,\lambda_n)=C_{n,\alpha,\eta} \prod_{i=1}^n\lambda_i^\alpha e^{-\lambda_i} \Delta_n^2(\boldsymbol{\lambda}) \sum_{k=1}^n \frac{e^{c_\eta \lambda_k}}{\displaystyle \prod_{\substack{i=1\\i\neq k}} ^n(\lambda_k-\lambda_i)}
  \end{align}
  where
  \begin{align*}
    C_{n,\alpha,\eta}=\frac{K_{n,\alpha}(n-1)!}{(\eta+1)^{\alpha+1}\eta^{n-1}},
  \end{align*}
  and $c_\eta=\eta/(\eta+1)$.
\end{corollary}


 {\color{blue}It is noteworthy that in addition to the contour integral approaches due to \cite{ref:wang,ref:peterBeta}}, the repeated application of the l'Hospital's rule due to \cite{ref:khatri}\footnote{Repeated
application of the l'Hospital's rule in the context of simplifying indeterminate forms involving determinants is given in \cite{ref:khatri}.} can be used to obtain the above form. Here we adopt the former approach, since it seems algebraically less tedious. As such, by substituting $\boldsymbol{\Sigma}=\textbf{I}_n+\eta \textbf{uu}^{*}$ into (\ref{eq_joint_pdf_ordered_eigenvals}) and simplifying\footnote{Here we make use of the relation $\left(\mathbf{I}_n+\eta \mathbf{uu}^*\right)^{-1}=\mathbf{I}_n-c_\eta \mathbf{uu}^*$.} the resultant joint p.d.f. with the help of (\ref{uniint}), we get
\begin{align}
\label{eq_joint_pdf_with_F}
    f(\lambda_1,\lambda_2,....,\lambda_n) &= \frac{K_{n,\alpha}}{(\eta+1)^{m}} \prod_{i=1}^{n}\lambda_i^\alpha e^{-\lambda_i} \Delta_n^2(\boldsymbol{\lambda}) \; _0\Tilde{F}_0\left(\boldsymbol{\Lambda},\mathbf{B}\right)
\end{align}{}
where $\mathbf{B}=\text{diag}\left(\frac{\eta}{\eta+1},0,...,0\right)$ is a rank-one matrix.
{\color{blue}Following the developments in \cite{ref:wang,ref:peterBeta}}, the object  $_0\Tilde{F}_0\left(\boldsymbol{\Lambda},\mathbf{B}\right)$ can be further simplified to yield
\begin{align}
    _0\Tilde{F}_0\left(\boldsymbol{\Lambda},\mathbf{B}\right) &= (n-1)!\left(\frac{\eta+1}{\eta}\right)^{n-1} \sum_{i=1}^n \frac{e^{c_\eta \lambda_i}}{\displaystyle \prod_{\substack{j=1 \\ j\neq i}}^n\left(\lambda_i-\lambda_j\right)},
\end{align}
which upon substituting into (\ref{eq_joint_pdf_with_F}) gives the desired result.

The functional form given in (\ref{cor joint}) facilitates the use of classical orthogonal polynomial approach due to Mehta \cite{ref:mehta} in our subsequent derivations.

\begin{definition}
  For $\rho>-1$, the generalized Laguerre polynomial of degree $M$, $L^{(\rho)}_M(z)$, is given by \cite{ref:szego}
  \begin{equation}
    \label{lagdef}
    L^{(\rho)}_M(z)=\frac{(\rho+1)_M}{M!}\sum_{j=0}^{M}\frac{(-M)_j}{(\rho+1)_j}\frac{z^j}{j!},
  \end{equation}
  with its $k$\textsuperscript{th} derivative satisfying
  \begin{align}
    \label{lagderi}
    \frac{{\rm d}^k}{{\rm d}z^k}L^{(\rho)}_M(z)=(-1)^kL^{(\rho+k)}_{M-k}(z),
  \end{align}
  where $(a)_j=a(a+1)\ldots(a+j-1)$ with  $(a)_0=1$ denotes the Pochhammer symbol. Moreover, we have the following contiguity relationships \cite[Eq. 8.971.4]{ref:gradshteyn}:
  \begin{align}
    zL_M^{(\rho)}(z)&=(M+\rho)L^{(\rho-1)}_M(z)-(M+1)L^{(\rho-1)}_{M+1}(z)\label{lagcont}\\
    L_M^{\rho-1}(z)&=L_M^\rho (z)-L_{M-1}^\rho (z)\label{lagcont1}.
  \end{align}
\end{definition}

Finally, we use the following compact notation to represent the
determinant of an $N\times N$ block matrix:
\begin{equation}
  \begin{split}
    \det\left[a_{i,j}\;\; b_{i,k-2}\right]_{\substack{i=1,2,\ldots,N\\
    j=1,2\\
    k=3,4,\ldots,N}}&=\left|\begin{array}{cccccc}
    a_{1,1} & a_{1,2}& b_{1,1}& b_{1,2}& \ldots & b_{1,N-2}\\
      a_{2,1} & a_{2,2}& b_{2,1}& b_{2,2}& \ldots & b_{2,N-2}\\
      \vdots & \vdots & \vdots & \vdots &\ddots & \vdots \\
      a_{N,1} & a_{N,2}& b_{N,1}& b_{N,2}& \ldots & b_{N,N-2}
    \end{array}\right|.
  \end{split}
\end{equation}

\section{ Probability Density Function of $\kappa_{\text{SC}}^2(\mathbf{X})$}

Here we derive a closed form p.d.f. expression for $\kappa_{\text{SC}}^2(\mathbf{X})$ by leveraging the orthogonal polynomial techniques due to Mehta \cite{ref:mehta}.
By definition, the m.g.f. of $\kappa_{\text{SC}}^2(\mathbf{X})$ is
\begin{align*}
  \mathcal{M}_{\kappa_{\text{SC}}^2(\mathbf{X})}(s)=e^{-s}\int_{\mathcal{R}} e^{-s\frac{\sum_{j=2}^n\lambda_j}
  {\lambda_1}}
  f(\lambda_1,\ldots,\lambda_n)
  {\rm d}\lambda_1\cdots {\rm d}\lambda_n
\end{align*}
where $\mathcal{R}=\{0\leq \lambda_1\leq\cdots\leq \lambda_n<\infty\}$. For convenience, let us introduce the substitution $\lambda_1=x$ and rewrite the above multiple integral, keeping the integration with respect to $x$ last, as
\begin{align}
  \label{int sep}
  \mathcal{M}_{\kappa_{\text{SC}}^2(\mathbf{X})}(s)=e^{-s}\int_0^\infty \int_{\mathcal{R}_x} e^{-s\frac{\sum_{j=2}^n\lambda_j}
  {x}}
  f(x,\lambda_2,\ldots,\lambda_n)
  {\rm d}\lambda_2\cdots {\rm d}\lambda_n{\rm d}x
\end{align}
where $\mathcal{R}_x=\{x\leq\cdots\leq \lambda_n<\infty\}$. To facilitate further analysis, we may use the decomposition, $\Delta_n^2(\boldsymbol{\lambda})=\prod_{i=2}^n(x-\lambda_i)^2
\Delta_{n-1}^2(\boldsymbol{\lambda})$, to rewrite (\ref{cor joint}) as
\begin{multline}
  \label{res den}
  f(x,\lambda_2,\ldots,\lambda_n)=C_{n,\alpha,\eta}
  x^\alpha e^{-x}
  \prod_{i=2}^n\lambda_i^{\alpha}e^{-\lambda_i}(x-\lambda_i)^2
  \Delta_{n-1}^2(\boldsymbol{\lambda}) \\
  \times 
  \left(\frac{e^{c_\eta x}}{\displaystyle\prod_{i=2}^n\left(x-\lambda_i\right)}+\sum_{k=2}^n\frac{e^{c_\eta \lambda_k}}{(\lambda_k-x)\displaystyle\prod_{\substack{i=2\\i\neq k}}^n\left(\lambda_k-\lambda_i\right)}\right).
\end{multline}
Therefore, we use (\ref{res den}) in (\ref{int sep}) with some algebraic manipulation to yield
\begin{align}
  \label{mgf decom}
  \mathcal{M}_{\kappa_{\text{SC}}^2(\mathbf{X})}(s) = \mathcal{A}(s) + \mathcal{B}(s)
  \end{align}
where
\begin{align}
  \label{Pdef}
  \mathcal{A}(s)=C_{n,\alpha,\eta} e^{-s}
  \int_0^\infty 
  e^{-x(1-c_\eta)}x^\alpha \Biggl( \int_{\mathcal{R}_x}
  \prod_{i=2}^n e^{-\left(1+\frac{s}{x}\right)\lambda_i}\lambda_i^\alpha (x-\lambda_i) 
  \Delta_{n-1}^2(\boldsymbol{\lambda})
  {\rm d}\lambda_2\cdots{\rm d}\lambda_n \Biggr)\; {\rm d}x
\end{align}
and 
\begin{multline}
  \label{Sdef}
  \mathcal{B}(s)=C_{n,\alpha,\eta} e^{-s} 
  \int_0^\infty e^{-x}x^{\alpha} \Biggl( \int_{\mathcal{R}_x}
  \sum_{k=2}^n\frac{e^{c_\eta \lambda_k}}{(\lambda_k-x)\displaystyle\prod_{\substack{i=2\\i\neq k}}^n\left(\lambda_k-\lambda_i\right)} \\
  \times \prod_{i=2}^ne^{-\left(1+\frac{s}{x}\right)\lambda_i}\lambda_i^{\alpha}(x-\lambda_i)^2
  \Delta_{n-1}^2(\boldsymbol{\lambda}){\rm d}\lambda_2\cdots{\rm d}\lambda_n\Biggr) {\rm d}x.
\end{multline}
Since the above two multiple integrals are structurally different from each other, in what follows, we evaluate each separately.

First, let us focus on $\mathcal{A}(s)$. As such, noting that the inner $(n-1)$-fold integral is symmetric in the variables $\lambda_2,\lambda_3,\ldots,\lambda_n$, we may remove the ordered region of integration to obtain
\begin{multline}
  \mathcal{A}(s) = \frac{C_{n,\alpha,\eta}}{(n-1)!}
  e^{-s} \int_0^\infty
  e^{-x(1-c_\eta)}x^\alpha \Bigg( \int_{[x,\infty)^{n-1}}
  \prod_{i=2}^n e^{-\left(1+\frac{s}{x}\right)\lambda_i}\lambda_i^\alpha (x-\lambda_i)\\ 
  \times\Delta_{n-1}^2(\boldsymbol{\lambda})
  {\rm d}\lambda_2\cdots{\rm d}\lambda_n\Bigg)\; {\rm d}x.
\end{multline}
To facilitate further analysis, we apply the change of variables, $y_{i-1}=(x+s)(\lambda_i-x)/x,\; i=2,\ldots,n$, to the inner $(n-1)$-fold integral with some algebraic manipulation to yield
\begin{multline}
    \mathcal{A}(s)=(-1)^{(n-1)(1+\alpha)}\frac{C_{n,\alpha,\eta}}{(n-1)!} e^{-ns}\int_0^\infty e^{-x(n-c_\eta)}\frac{x^{n(n+\alpha-1)}}{(x+s)^{(n-1)(n+\alpha)}} 
    R^{(\alpha)}_{n-1}\left(-(x+s)\right){\rm d}x
\end{multline}
where
\begin{align*}
  R^{(\alpha)}_n(t)=\int_{[0,\infty)^n}\prod_{j=1}^n e^{-y_j}y_j(t-y_j)^\alpha \Delta^2_n(\mathbf{y}){\rm d}y_1\cdots{\rm d}y_n.
\end{align*}
The above integral can be solved using the powerful orthogonal polynomial technique devised in \cite[Section 22.2.2]{ref:mehta} to yield (see also \cite[Eq. c.6]{ref:prathapJMVA})
\begin{align}
  \label{eq_R_eval}
  R^{(\alpha)}_n(t)=(-1)^{n\alpha}\prod_{j=0}^{n-1}(j+1)!(j+1)!\prod_{j=0}^{\alpha-1}\frac{(n+j)!}{j!}
  \det\left[L^{(j)}_{n+i-j}(t)\right]_{i,j=1,\ldots,\alpha}.
\end{align}
This in turn enables us to write $\mathcal{A}(s)$ as
\begin{multline}
  \label{eq_Afinal}
  \mathcal{A}(s)=\tilde{K}_{n,\alpha,\eta} e^{-ns}\int_0^\infty e^{-x(n-c_\eta)}\frac{x^{n(n+\alpha-1)}}{(x+s)^{(n-1)(n+\alpha)}} \det\left[L_{n+i-j-1}^{(j)}(-x-s) \;\right]_{i,j=1,..,\alpha} \; {\rm d}x
\end{multline}
where $\tilde{K}_{n,\alpha,\eta}=(-1)^{n-1}\frac{(n-1)!}{(n+\alpha-1)!(\eta+1)^{\alpha+1}\eta^{n-1}}$. For convenience, let us leave $\mathcal{A}(s)$ without further simplification and instead focus on $\mathcal{B}(s)$.

Due to symmetry, we can convert the ordered region of integration into an unordered region to yield
\begin{multline}
  \label{Bman}
  \mathcal{B}(s)=\frac{C_{n,\alpha,\eta}}{(n-1)!} e^{-s}
  \int_0^\infty e^{-x}x^{\alpha} \Biggl(\int_{[x,\infty)^{n-1}}
  \sum_{k=2}^n\frac{e^{c_\eta \lambda_k}}{(\lambda_k-x)\displaystyle\prod_{\substack{i=2\\i\neq k}}^n\left(\lambda_k-\lambda_i\right)} \\
  \times \prod_{i=2}^ne^{-\left(1+\frac{s}{x}\right)\lambda_i}\lambda_i^{\alpha}(x-\lambda_i)^2 
  \Delta_{n-1}^2(\boldsymbol{\lambda}){\rm d}\lambda_2\cdots{\rm d}\lambda_n\Bigg) {\rm d}x.
\end{multline}
A careful inspection of the summation in the inner integrand reveals that each term therein contributes the same amount to the final answer. Therefore, capitalizing on that observation, we can further simplify the above multiple integral to obtain
\begin{multline}
  \label{Bman1}
  \mathcal{B}(s)=\frac{C_{n,\alpha,\eta}}{(n-2)!}
  e^{-s}
  \int_0^\infty e^{-x}x^{\alpha} \Biggl(\int_{[x,\infty)^{n-1}}
  \frac{e^{c_\eta \lambda_2}}{(\lambda_2-x)\displaystyle\prod_{i=3}^n\left(\lambda_2-\lambda_i\right)} \\
  \times \prod_{i=2}^ne^{-\left(1+\frac{s}{x}\right)\lambda_i}\lambda_i^{\alpha}(x-\lambda_i)^2
  \Delta_{n-1}^2(\boldsymbol{\lambda}){\rm d}\lambda_2\cdots{\rm d}\lambda_n\Bigg) {\rm d}x.
\end{multline}
Noting the decomposition, $\Delta_{n-1}^2(\boldsymbol{\lambda})=\prod_{j=3}^n(\lambda_2-\lambda_j)^2\Delta_{n-2}^2(\boldsymbol{\lambda})$, the above multiple integral can be rewritten as
\begin{multline}
  \label{Bman2}
  \mathcal{B}(s)=\frac{C_{n,\alpha,\eta}}{(n-2)!}
  e^{-s}
  \int_0^\infty e^{-x}x^{\alpha} \Biggl\{\int_{x}^\infty 
  e^{-\left(1+\frac{s}{x}-c_\eta\right)\lambda_2}\lambda_2^\alpha (\lambda_2-x)\\
  \times \Biggl( \int_{[x,\infty)^{n-2}} \prod_{i=3}^ne^{-\left(1+\frac{s}{x}\right)\lambda_i}\lambda_i^{\alpha} (x-\lambda_i)^2(\lambda_2-\lambda_i)
  \Delta_{n-2}^2(\boldsymbol{\lambda}){\rm d}\lambda_3\cdots{\rm d}\lambda_n\Biggr) {\rm d}\lambda_2\Biggr\}{\rm d}x.
\end{multline}
Now it is convenient to introduce the variable transformations, $y=\lambda_2-x$ and $y_{i-2}=(x+s)(\lambda_i-x)/x$, $i=3,\ldots,n$, in the above multiple integral to yield 
\begin{multline}
  \label{Bman3}
  \mathcal{B}(s)=(-1)^{n\alpha}\frac{C_{n,\alpha,\eta}}{(n-2)!}
  e^{-ns}
  \int_0^\infty \frac{e^{-x(n-c_\eta)}x^{\alpha}}{\left(1+\frac{s}{x}\right)^{(n-2)(n+\alpha+1)}} \Biggl\{\int_{0}^\infty 
  e^{-\left(1+\frac{s}{x}-c_\eta\right)y} y (y+x)^\alpha \\
  \times T^{(\alpha)}_{n-2}\left(\left(1+\frac{s}{x}\right)y,-s-x\right){\rm d}y\Biggr\} {\rm d}x
\end{multline}
where
\begin{align}
    T^{(\alpha)}_n(a,b):=\int_{[0,\infty)^n}\prod_{i=1}^n(a-y_i)(b-y_i)^\alpha e^{-y_i}y_i^2 \Delta_n^2(\mathbf{y}) {\rm d}y_1\cdots{\rm d}y_n.
\end{align}
Following \cite[Eq. 5.8]{ref:prathapSIAM}, the above multiple integral admits the solution
\begin{align}
  \label{Teval}
  T^{(\alpha)}_n(a,b):=\frac{(-1)^{n+\alpha(n+\alpha)}\overline{\mathcal{K}}_{n,\alpha}}{(b-a)^\alpha}\det\left[L^{(2)}_{n+i-1}(a)\;\;\; L_{n+i+1-j}^{(j)}(b)\right]_{\substack{i=1,\ldots,\alpha+1\\j=2,\ldots,\alpha+1}}
\end{align}
where 
\begin{align*}
  \overline{\mathcal{K}}_{n,\alpha}=\frac{\prod_{j=1}^{\alpha+1}(n+j-1)! \prod_{j=0}^{n-1}(j+1)!(j+2)!}{\prod_{j=0}^{\alpha-1}j!}.
\end{align*}
Therefore, we use (\ref{Teval}) in (\ref{Bman3}) with some algebraic manipulation to arrive at
\begin{multline*}
  \mathcal{B}(s)
  =-\tilde{K}_{n,\alpha,\eta} e^{-ns} \int_0^\infty \frac{e^{-x(n-c_\eta)}x^{\alpha}}{\left(1+\frac{s}{x}\right)^{(n-1)(n+\alpha)-2}} \left\{\int_0^\infty y  e^{-\left(1+\frac{s}{x}-c_\eta\right)y}
  \right.\\
  \times\det\left[L^{(2)}_{n+i-3}\left(y\left(1+\frac{s}{x}\right)\right)\;\;\; L_{n+i-1-j}^{(j)}\left(-x-s\right)\right]_{\substack{i=1,\ldots,\alpha+1\\j=2,\ldots,\alpha+1}} {\rm d}y\Biggr\}
  {\rm d}x
\end{multline*}
from which we obtain, after the variable transformation $y(1+s/x)=t$, 
\begin{multline*}
  \mathcal{B}(s)
  =-\tilde{K}_{n,\alpha,\eta} e^{-ns} \int_0^\infty \frac{e^{-x(n-c_\eta)}x^{\alpha}}{\left(1+\frac{s}{x}\right)^{(n-1)(n+\alpha)}} \left\{\int_0^\infty t  e^{-t\left(1-\frac{c_\eta x}{x+s}\right)}
  \right.\\
  \times\det\left[L^{(2)}_{n+i-3}\left(t\right)\;\;\; L_{n+i-1-j}^{(j)}\left(-x-s\right)\right]_{\substack{i=1,\ldots,\alpha+1\\j=2,\ldots,\alpha+1}} {\rm d}t\Biggr\}
  {\rm d}x.
\end{multline*}
Since only the first column of the determinant depends on $t$, we can easily rewrite the above double integral as
\begin{multline}
  \label{eq_B_tint}
  \mathcal{B}(s)
  =-\tilde{K}_{n,\alpha,\eta} e^{-ns}\int_0^\infty e^{-x(n-c_\eta)}\frac{x^{n(n+\alpha-1)}}{(x+s)^{(n-1)(n+\alpha)}} \\
  \times \det[\rho_i(x,s) \;\; L_{n+i-j-1}^{(j)}(-x-s) \;]_{\substack{i=1,\ldots,\alpha+1\\j=2,\ldots,\alpha+1}} \; {\rm d}x
\end{multline}
where
\begin{align}
  \label{rhodef}
  \rho_i(x,s) &=  \int_0^\infty t  e^{-t\left(1-\frac{c_\eta x}{x+s}\right)} L_{n+i-3}^{(2)}(t) {\rm d}t.
  \end{align}
Following (\ref{lagcont}), we further decompose $\rho_i(s,x)$ to yield
\begin{multline}
  \rho_i(x,s) = (n+i-1)\int_0^\infty e^{-t\left(1-\frac{c_\eta x}{x+s}\right)} L_{n+i-3}^{(1)}(t){\rm d}t   - (n+i-2)\int_0^\infty e^{-t\left(1-\frac{c_\eta x}{x+s}\right)} L_{n+i-2}^{(1)}(t){\rm d}t,
\end{multline}
from which we obtain in view of \cite[Eq. 7.414.5]{ref:gradshteyn}
\begin{multline}
  \rho_i(x,s) = (n+i-1)\sum_{m=0}^{n+i-3}\frac{\left(-\frac{c_\eta x}{x+s}\right)^{n+i-3-m}}{\left(1-\frac{c_\eta x}{x+s}\right)^{n+i-2-m}}   - (n+i-2)\sum_{m=0}^{n+i-2}\frac{\left(-\frac{c_\eta x}{x+s}\right)^{n+i-2-m}}{\left(1-\frac{c_\eta x}{x+s}\right)^{n+i-1-m}}.
\end{multline}
Some algebraic manipulation now gives
\begin{align}
  \label{eq_rho_decom}
  \rho_i(x,s) &= 1 + \sigma_i(x+s,x)
\end{align}
where
\begin{align*}
  \sigma_i(x+s,x) &= (-1)^{n+i-3}\left(\frac{c_\eta x}{x+s-c_\eta x}\right)^{n+i-2}\left(\frac{(n+i-1)(x+s)-c_\eta x}{x+s-c_\eta x}\right).
\end{align*} 
Now it is convenient to substitute (\ref{eq_rho_decom}) into (\ref{eq_B_tint}) and exploit the multilinear property to expand the resultant determinant to yield
\begin{multline}
\label{eq_B_rhosub}
  \mathcal{B}(s) = -\tilde{K}_{n,\alpha,\eta} e^{-ns}\int_0^\infty e^{-x(n-c_\eta)}\frac{x^{n(n+\alpha-1)}}{(x+s)^{(n-1)(n+\alpha)}} \det[1 \;\; L_{n+i-j-1}^{(j)}(-x-s) \;]_{\substack{i=1,\ldots,\alpha+1\\j=2,\ldots,\alpha+1}} \; {\rm d}x \\
  -\tilde{K}_{n,\alpha,\eta} e^{-ns}\int_0^\infty e^{-x(n-c_\eta)}\frac{x^{n(n+\alpha-1)}}{(x+s)^{(n-1)(n+\alpha)}} \\
  \times \det[\sigma_i(x+s,x) \;\; L_{n+i-j-1}^{(j)}(-x-s) \;]_{\substack{i=1,\ldots,\alpha+1\\j=2,\ldots,\alpha+1}} \; {\rm d}x.
\end{multline}
In order to further simplify the determinant in the first integral, we apply the following row operations
\begin{align*}
  i\text{th row}\to i\text{th row}+(-1)(i-1)\text{th row},\;\; i=2,3,\ldots,\alpha+1
\end{align*}
and expand the resultant determinant using its first column to obtain
\begin{align}
\label{eq Lag det simp}
  \det\left[ 1\;\;\;L_{n+i-1-j}^{(j)}\left(-x-s\right)\right]_{\substack{i=1,\ldots,\alpha+1\\j=2,\ldots,\alpha+1}}&
  =\det\left[ L_{n+i-1-j}^{(j)}\left(-x-s\right)\right]_{i,j=1,\ldots,\alpha}
\end{align}
where we have made use of the contiguous relation (\ref{lagcont1}). Therefore, in view of (\ref{eq_Afinal}), (\ref{eq_B_rhosub}) can be written as
\begin{multline}
  \label{eq_B_detexp}
  \mathcal{B}(s)=
  -\mathcal{A}(s)
  -\tilde{K}_{n,\alpha,\eta} e^{-ns} \int_0^\infty e^{-x(n-c_\eta)}\frac{x^{n(n+\alpha-1)}}{(x+s)^{(n-1)(n+\alpha)}} \\
  \times 
  \det[\sigma_i(x+s,x) \;\; L_{n+i-j-1}^{(j)}(-x-s) \;]_{\substack{i=1,\ldots,\alpha+1\\j=2,\ldots,\alpha+1}} \; {\rm d}x,
\end{multline}
from which we obtain, following (\ref{mgf decom})
\begin{multline}
  \label{eq_mgf_final}
  \mathcal{M}_{\kappa_{\text{SC}}^2(\mathbf{X})}(s)=
  -\tilde{K}_{n,\alpha,\eta} e^{-ns} \int_0^\infty e^{-x(n-c_\eta)}\frac{x^{n(n+\alpha-1)}}{(x+s)^{(n-1)(n+\alpha)}}\\
  \times 
  \det[\sigma_i(x+s,x) \;\; L_{n+i-j-1}^{(j)}(-x-s) \;]_{\substack{i=1,\ldots,\alpha+1\\j=2,\ldots,\alpha+1}} \; {\rm d}x.
\end{multline}
Keeping in mind that $\mathcal{L}^{-1}\left\{e^{-ns}\,G(x+s)\right\}=e^{nx}e^{-xz}\mathcal{L}^{-1}\left\{e^{-ns}\,G(s)\right\}$, we take the inverse Laplace transform of (\ref{eq_mgf_final}) to obtain
\begin{multline}
  \label{eq_mgf_inv}
  f^{\alpha}_{\kappa_{\text{SC}}^2(\mathbf{X})}(z)=
  -\tilde{K}_{n,\alpha,\eta} \mathcal{L}^{-1}\left\{\frac{e^{-ns}}{s^{(n-1)(n+\alpha)}}\int_0^\infty e^{-x(z-c_\eta)}x^{n(n+\alpha-1)}\right.\\
  \left. \times \det\left[\sigma_i(s,x) \;\; L_{n+i-j-1}^{(j)}(-s) \;\right]_{\substack{i=1,\ldots,\alpha+1\\j=2,\ldots,\alpha+1}} \; {\rm d}x\right\},
\end{multline}
from which, one obtains after changing the order of integration
\begin{align}
  \label{eq_fz_with_new}
  f_{\kappa_{\text{SC}}^2(\mathbf{X})}^\alpha(z) 
  &= -\tilde{K}_{n,\alpha,\eta} \; \int_0^\infty e^{-x(z-c_\eta)} x^{n(n+\alpha-1)}  \mathcal{L}^{-1}\left\{ D(x,s)\right\} {\rm d}x
\end{align}
where
\begin{align}
  D(x,s)=\det \left[A_i(x,s) \;\;\; L_{n+i-j-1}^{(j)}(-s)\right]_{\substack{i=1,\dotsc,\alpha+1 \\ j=2,\dotsc,\alpha+1}}
\end{align}
with
\begin{align}
  A_i(x,s) &= \frac{e^{-ns}\sigma_i(s,x)}{s^{(n-1)(n+\alpha)}}
\end{align}
and $\mathcal{L}^{-1}(\cdot)$ denotes the inverse Laplace transform.

Now let us focus on further simplification of $D(x,s)$. To this end, we use (\ref{lagdef})
to rewrite $D(x,s)$ as
\begin{multline}
\label{lagexpandeq}
  D(x,s) = \det \left[A_i(x,s) \;\;\;\;\;\; \frac{(j+1)_{n+i-j-1}}{(n+i-j-1)!} \sum_{k_j=0}^{n+i-j-1}\frac{(-n-i+j+1)_{k_j}(-s)^{k_j}}{(j+1)_{k_j}k_j!} \right]_{\substack{i=1,..,\alpha+1 \\ j=2,..,\alpha+1}}.
\end{multline}
Further simplification of the determinant is difficult in its current form due to the dependence of the upper limit of the finite summation on $i$ and $j$. To circumvent this challenge, 
we use the decomposition 
\begin{align*}
  (-n-i+j+1)_{k_j}&= (-n-i+j+1)_{k_j}\frac{(-n-\alpha+j)_{k_j}}{(-n-\alpha+j)_{k_j}} \\
  &= \frac{(n+i-j-1)!(n+\alpha-j-k_j)!(-n-\alpha+j)_{k_j}}{\Gamma(n+i-j-k_j)(n+\alpha-j)!}
\end{align*}
 with some algebraic manipulation to rewrite
\begin{align}
  & D(x,s)\nonumber \\
  &= \text{det} \Biggl[A_i(x,s) \qquad \frac{(n+i-1)!}{(n+\alpha-j)!j!} \nonumber\\
  &  \qquad \qquad \qquad \qquad \quad \;\; \times \sum_{k_j=0}^{n+\alpha-j}\frac{(-n-\alpha+j)_{k_j}(-s)^{k_j}}{(j+1)_{k_j}k_j!}\frac{(n+\alpha-j-k_j)!}{\Gamma(n+i-j-k_j)} \Biggr]_{\substack{i=1,..,\alpha+1 \\ j=2,..,\alpha+1}} \nonumber\\
  &= \frac{\prod_{j=1}^{\alpha+1}(n+j-1)!}{\prod_{j=2}^{\alpha+1} j!(n+\alpha-j)!} \sum_{k_2=0}^{n+\alpha-2} \text{...} \sum_{k_{\alpha+1}=0}^{n-1}\prod_{j=2}^{\alpha+1}\frac{(-n-\alpha+j)_{k_j}(-s)^{k_j}(n+\alpha-j-k_j)!}{(j+1)_{k_j}k_j!} \nonumber\\
  & \qquad \qquad \qquad \qquad \qquad \qquad \qquad \times \det \left[\frac{A_i(x,s)}{(n+i-1)!} \;\;\; \frac{1}{\Gamma(n+i-j-k_j)} \right]_{\substack{i=1,..,\alpha+1 \\ j=2,..,\alpha+1}}
  \end{align}
  where $\Gamma(\cdot)$ denotes the Gamma function.
  It is worth noting that the original restriction on $k_j$ in (\ref{lagexpandeq}), $k_j\leq n+i-j-1$, is implicitly imposed here by the term $\frac{1}{\Gamma(n+i-j-k_j)}$. Now we may collect all powers of $s$ in the outer nested summations and multiply the first column of the determinant by the resultant term $s^{\sum_{j=2}^{\alpha+1} k_j}$ to obtain
  \begin{align}
  \label{eq_D}
  & D(x,s)\nonumber \\
  &= \frac{(n+\alpha-1)!(n+\alpha)!}{(n-1)!} \sum_{k_2=0}^{n+\alpha-2} \text{...} \sum_{k_{\alpha+1}=0}^{n-1}\prod_{j=2}^{\alpha+1}\frac{(n+\alpha-j)!}{(j+k_j)!k_j!} \nonumber \\
  &\quad \qquad \qquad \qquad \qquad \qquad \times \det \left[\frac{A_i(x,s)s^{\sum_{j=2}^{\alpha+1}k_j}}{(n+i-1)!} \;\;\; \frac{1}{\Gamma(n+i-j-k_j)} \right]_{\substack{i=1,..,\alpha+1 \\ j=2,..,\alpha+1}}.
\end{align}
Again, the same restriction on the parameters $k_j$ is critically important for the existence of the term $s^{\sum_{j=2}^{\alpha+1} k_j}$. However, in what follows, for the clarity of presentation, we tacitly avoid it, since this particular restriction is implicitly embedded in the term $1/\Gamma(n+i-j-k_j)$ .   
 Noting that only the first column of the above determinant contains $s$, we take term-by-term Laplace inversion to obtain 
\begin{multline}
  \label{Lapinvdet}
  \mathcal{L}^{-1} \left\{D(x,s)\right\} = \frac{(n+\alpha-1)!(n+\alpha)!}{(n-1)!} \sum_{k_2=0}^{n+\alpha-2} \dotsc \sum_{k_{\alpha+1}=0}^{n-1} \prod_{j=2}^{\alpha+1} \frac{(n+\alpha-j)!}{(j+k_j)!k_j!} \\
  \times \det \left[B_i(x,z) \;\;\; \frac{1}{\Gamma(n+i-j-k_j)}\right]_{\substack{i=1,\dotsc,\alpha+1 \\ j=2,\dotsc,\alpha+1}}
\end{multline}
where 
\begin{align}
\label{Bdef}
    &B_i(x,z)\nonumber\\
    &= \mathcal{L}^{-1} \left\{\frac{A_i(x,s)s^{\sum_{j=2}^{\alpha+1}k_j}}{(n+i-1)!}\right\} \nonumber \\
    &= \frac{(-1)^{n+i-3}\left(c_\eta x\right)^{n+i-2}}{(n+i-1)!}\left((n+i-1) \mathcal{L}^{-1} \left\{ \frac{e^{-ns}}{s^{(n-1)(n+\alpha)-\sum_{j=2}^{\alpha+1}k_j-1}(s-c_\eta x)^{n+i-1}}\right\} \right. \nonumber \\
    &\qquad\qquad\qquad\qquad\qquad\qquad\qquad \left. - c_\eta x \mathcal{L}^{-1} \left\{ \frac{e^{-ns}}{s^{(n-1)(n+\alpha)-\sum_{j=2}^{\alpha+1}k_j}(s-c_\eta x)^{n+i-1}} \right\}\right).
\end{align}
Consequently, we make use of the Laplace inversion relation \cite[Eq. 6.10.6]{ref:erdelyi}
\begin{align}
  \mathcal{L}^{-1} \left\{\frac{e^{-ns}}{s^a (s-\omega)^b} \right\} &= \frac{(z-n)^{a+b-1}}{\Gamma(a+b)} \;_1F_1\left(b;a+b;\omega(z-n)\right)H(z-n),\;\; a+b>0,
\end{align}
with ${}_1F_1 (\cdot;\cdot;\cdot)$ denoting the confluent hypergeometric function of the first kind \cite{ref:gradshteyn} and $H(\cdot)$ denoting the Heaviside unit step function, in (\ref{Bdef}) with some algebraic manipulation to obtain
\begin{multline}
  \label{Bdefalt}
  B_i(x,z) = \frac{(-1)^{a_i}(c_\eta x)^{a_i-1}(z-n)^{b_i-1}}{\Gamma(a_i)\Gamma(b_i)} \Biggl( \;_1F_1\left(a_i;b_i;\mathcal{Z}_x\right)-\frac{\mathcal{Z}_x}{a_i b_i}\;_1F_1\left(a_i;b_i+1;\mathcal{Z}_x\right) \Biggr) H(z-n)
\end{multline}
where $a_i = n+i-1$, $b_i = n^2+n\alpha+i-\alpha-\sum_{j=2}^{\alpha+1}k_j-2$, and $\mathcal{Z}_x = c_\eta x(z-n)$. We note here that $b_i>0$ by virtue of $k_j\leq n-i-j-1$. Now it is convenient to substitute (\ref{Bdefalt}) into (\ref{Lapinvdet}) and use the resultant expression in (\ref{eq_fz_with_new}) to yield
\begin{multline}
  f_{\kappa_{\text{SC}}^2(\mathbf{X})}^\alpha(z) = -\tilde{K}_{n,\alpha,\eta} \frac{(n+\alpha-1)!(n+\alpha)!}{(n-1)!} \sum_{k_2=0}^{n+\alpha-2} \dotsc \sum_{k_{\alpha+1}=0}^{n-1} \prod_{j=2}^{\alpha+1} \frac{(n+\alpha-j)!}{(j+k_j)!k_j!} \\
  \times \int_0^\infty e^{-x(z-c_\eta)} x^{n(n+\alpha-1)}
  \det \left[B_i(x,z) \;\;\; \frac{1}{\Gamma(n+i-j-k_j)}\right]_{\substack{i=1,\dotsc,\alpha+1 \\ j=2,\dotsc,\alpha+1}} {\rm d}x.
\end{multline}
Since only the first column of the determinant depends on $x$, the integration operation can be absorbed into the determinant to obtain
\begin{multline}
    \label{eq_fz_with_Ii_method2}
    f_{\kappa_{\text{SC}}^2(\mathbf{X})}^\alpha(z) = -\tilde{K}_{n,\alpha,\eta} \frac{(n+\alpha-1)!(n+\alpha)!}{(n-1)!} \sum_{k_2=0}^{n+\alpha-2} \dotsc \sum_{k_{\alpha+1}=0}^{n-1} \prod_{j=2}^{\alpha+1} \frac{(n+\alpha-j)!}{(j+k_j)!k_j!} \\
    \times  \det \left[\mathcal{I}_i(z) \;\;\; \frac{1}{\Gamma(n+i-j-k_j)}\right]_{\substack{i=1,\dotsc,\alpha+1 \\ j=2,\dotsc,\alpha+1}}
\end{multline}
where
\begin{align}
  \mathcal{I}_i(z) &= \int_0^\infty e^{-x(z-c_\eta)} x^{n(n+\alpha-1)} B_i(x,z) {\rm d}x \nonumber \\
  &= \frac{(-1)^{a_i}(c_\eta)^{a_i-1}(z-n)^{b_i-1}}{\Gamma(a_i)\Gamma(b_i)}\left( \int_0^\infty e^{-x(z-c_\eta)} x^{c_i-1}\;_1F_1\left(a_i;b_i;\mathcal{Z}_x\right) {\rm d}x\right. \nonumber \\
  &\qquad\qquad\qquad\quad \left. - \frac{c_\eta(z-n)}{a_i b_i} \int_0^\infty e^{-x(z-c_\eta)} x^{c_i}\;_1F_1\left(a_i;b_i+1;\mathcal{Z}_x\right) {\rm d}x \right) H(z-n)
\end{align}
with $c_i=n^2+n\alpha+i-1$.
Now in light of \cite[Eq. 7.621.4]{ref:gradshteyn} each of the above integrals can be evaluated in closed-form to obtain
\begin{multline}
  \label{eq I}
  \mathcal{I}_i(z) = \frac{(-1)^{a_i}\Gamma(c_i) c_\eta^{a_i-1}(z-n)^{b_i-1}}{\Gamma(a_i)\Gamma(b_i)(z-c_\eta)^{c_i}} \left(\;_2F_1\left(a_i,c_i;b_i;\frac{c_\eta(z-n)}{z-c_\eta}\right)\right. \\
  \left. - \frac{c_i}{a_ib_i}\frac{c_\eta(z-n)}{z-c_\eta}\;_2F_1\left(a_i,c_i+1;b_i+1;\frac{c_\eta(z-n)}{z-c_\eta}\right)\right) H(z-n)
\end{multline}
where ${}_2F_1(a,b;c;z)$ denotes the Gauss hypergeometric function \cite{ref:gradshteyn}. 
In order to further simplify the above expression, noting that ${}_2F_1(a,b;c;z)={}_3F_2(a,b,d;c,d;z)$ for $d\neq 0$, we re-write (\ref{eq I}) as
\begin{align}
  \mathcal{I}_i & (z) \nonumber \\ = & \frac{(-1)^{a_i}\Gamma(c_i) c_\eta^{a_i-1}(z-n)^{b_i-1}}{\Gamma(a_i)\Gamma(b_i)(z-c_\eta)^{c_i}} H(z-n) \Biggl(\;_3F_2\left(a_i,c_i,a_i+1;a_i+1,b_i;\frac{c_\eta(z-n)}{z-c_\eta}\right) \nonumber \\
   &- \frac{(a_i-(a_i-1))c_i (a_i+1)}{(a_i+1)a_ib_i}\frac{c_\eta(z-n)}{z-c_\eta} {}_3F_2\left(a_i,c_i+1,a_i+2;a_i+2,b_i+1;\frac{c_\eta(z-n)}{z-c_\eta}\right)\Biggr) \text{.}
\end{align}
Now, the following contiguous relationship \cite[Eq. 07.27.17.0015.01]{wolfram}
\begin{multline}
  {}_3F_2\left(a,b,c;d,e;z\right) = {}_3F_2\left(a+1,b,c;d+1,e;z\right) \\
  - \frac{(d-a)bcz}{(d+1)de}{}_3F_2\left(a+1,b+1,c+1;d+2,e+1;z\right)
\end{multline}
can be used to further simplify the above expression to arrive at
\begin{multline}
  \mathcal{I}_i(z) = \frac{(-1)^{a_i}\Gamma(c_i)c_\eta^{a_i-1}(z-n)^{b_i-1}}{\Gamma(a_i)\Gamma(b_i)(z-c_\eta)^{c_i}}  {}_3F_2\left(a_i-1,c_i,a_i+1;a_i,b_i;\frac{c_\eta(z-n)}{z-c_\eta}\right) H(z-n).
\end{multline}
Substituting this result into (\ref{eq_fz_with_Ii_method2}), along with some algebraic manipulations and the index shift $i\rightarrow i-1, j\rightarrow j-1$ gives us the final result which is given by the following theorem.

\begin{theorem}
  \label{thm exact pdf}
  The exact p.d.f. of $\kappa_{\text{SC}}^2(\mathbf{X})$  is given by
  \begin{multline}
    \label{eq_pdf_Z}
    f_{\kappa_{\text{SC}}^2(\mathbf{X})}^\alpha(z) = \frac{(n+\alpha)!(z-n)^{n^2+n\alpha-\alpha-2}}{(\eta+1)^{n+\alpha}(z-c_\eta)^{n^2+n\alpha}} \sum_{k_1=0}^{n+\alpha-2} \dotsc  \sum_{k_{\alpha}=0}^{n-1}\left(\prod_{j=1}^{\alpha}\frac{(n+\alpha-j-1)!}{(j+k_j+1)!k_j!} \right)\\
    \times (z-n)^{-\sum_{j=1}^{\alpha}k_j}
    \mathrm{det}\left[ \mathcal{G}_i(z,\eta) \;\; \frac{1}{\Gamma(n+i-j-k_j)} \right]_{\substack{i=0,..,\alpha \\ j=1,..,\alpha}}H(z-n)
  \end{multline}
  where
  \begin{multline}
  \label{g_i_in_hypergeometric}
    \mathcal{G}_i(z,\eta) = (-1)^i\frac{c_\eta^i(z-n)^i}{(z-c_\eta)^i} \frac{\Gamma\left(n^2+n\alpha+i\right)}{\Gamma(n+i)\Gamma\left(n^2+n\alpha-\alpha+i-\sum_{j=1}^{\alpha}k_j-1\right)}  \\ 
    \times {}_3F_2 \Biggl( n+i-1,n+i+1,n^2+n\alpha+i; n+i, n^2+n\alpha-\alpha+i-\sum_{j=1}^{\alpha}k_j-1; \frac{c_\eta(z-n)}{z-c_\eta} \Biggr)
  \end{multline}{}
  and $H(z)$ is the unit step function.
\end{theorem}
\begin{remark}
{\color{blue} It is worth mentioning that the  above generalized hypergeometric function ${}_3F_2$
boils down to a sum of rational functions, thereby simplifying the classical equivalent infinite series expansion. To demonstrate this, let us 
utilize the decomposition\footnote{Capitalizing on the observation $\frac{(a)_k}{(a-1)_k}=1+\frac{k}{(a-1)}$, we obtain ${}_3F_2(a,b,c;a-1,d;z)= \displaystyle\sum_{k=0}^\infty \left(1+\frac{k}{a-1}\right) \frac{(b)_k (c)_k}{(d)_k k!}z^k={}_2F_1(b,c;d;z)+\frac{z}{a-1}\displaystyle\sum_{k=0}^\infty \frac{(b)_{k+1} (c)_{k+1}}{(d)_{k+1}k!)!}z^k$. Now the final result follows by  noting the decomposition $(p)_{k+1}=p(p+1)_k$.}
\begin{align}
    {}_3F_2(a,b,c;a-1,d;z)={}_2F_1(b,c;d;z)+\frac{zbc}{d(a-1)} {}_2F_1(b+1,c+1;d+1;z),
\end{align} along with \cite[Eq. 9.131.1]{ref:gradshteyn} to  further simplify the generalized hypergeometric function in (\ref{g_i_in_hypergeometric}) as
\begin{align}
\label{3F2decom}
    {}_3F_2&\left(p_i+1,p_i-1,q_i;p_i,r_i;\frac{c_\eta(z-n)}{z-c_\eta}\right)\nonumber \\
    &= \left(\frac{z(\eta+1)-\eta}{z+\eta(n-1)}\right)^{p_i-1} (q_i-r_i)! \sum_{k=0}^{q_i-r_i} \frac{(p_i-1)_k}{(q_i-r_i-k)!k!(r_i)_k} \left(\frac{\eta(z-n)}{z+\eta(n-1)}\right)^k  \nonumber \\
    & \hspace{0.4\textwidth} \times \left(1+\frac{(p_i-1+k)q_i}{(r_i+k)p_i}\frac{\eta(z-n)}{z+\eta(n-1)}\right)
\end{align}
where $p_i=n+i, q_i=n^2+n\alpha+i$, and $r_i=n^2+n\alpha-\alpha+i-\sum_{j=1}^\alpha k_j-1$. Nevertheless, for notational concision, we use ${}_3F_2$ instead of the above representation involving the sum of rational functions.
}
\end{remark}
  Since the number of nested summations depends only on $\alpha$, this formula provides an efficient way of evaluating the p.d.f. of $\kappa_{\text{SC}}^2(\mathbf{X})$, especially for small values of $\alpha$. Furthermore, since the algebraic complexity depends only on $n$ and $\alpha$ (i.e., the difference between $m$ and $n$), this in turn makes our result (\ref{eq_pdf_Z}) very useful for determining the macroscopic limit of $\kappa_{\text{SC}}^2(\mathbf{X})$ (i.e., as $m$ and $n$ grow large while their difference is fixed). This will be addressed in the next section.

Now it is worth observing that for some small values of $\alpha$, (\ref{eq_pdf_Z}) admits the following simple forms.

\begin{corollary}
  The exact p.d.f.s of $\kappa_{\text{SC}}^2(\mathbf{X})$ corresponding to $\alpha=0$ and $\alpha=1$ are given, respectively, by
  \begin{multline}
    \label{eq_pdf_alpha0}
    f_{\kappa_{\text{SC}}^2(\mathbf{X})}^0(z) = \frac{n(n^2-1)(z-n)^{n^2-2}}{(\eta+1)^n (z-c_\eta)^{n^2}}  {}_3F_2\left(n-1,n+1,n^2;n,n^2-1;\frac{c_\eta(z-n)}{z-c_\eta}\right) H(z-n),
  \end{multline}
  \begin{multline}
    f_{\kappa_{\text{SC}}^2(\mathbf{X})}^1(z) =\frac{(z-c_\eta)^{-n(n+1)}}{n(\eta+1)^{n+1}} \Biggl( \frac{n(n^2+n-1)!(3)_{n-1}}{(n-1)!} \sum_{j=0}^{n-1}\frac{(-n+1)_j(-1)^j(z-n)^{n^2+n-3-j}}{j!(3)_j(n^2+n-3-j)!}\\ 
    \left. \times_3F_2\left(n-1,n+1,n^2+n;n,n^2+n-2-j;\frac{c_\eta(z-n)}{z-c_\eta}\right)\right. \\
    \left.+ \left(\frac{c_\eta}{z-c_\eta}\right)\frac{(n+1)(n^2
    +n)!(3)_{n-2}}{(n-2)!}\sum_{j=0}^{n-2}\frac{(-n+2)_j(-1)^j(z-n)^{n^2+n-2-j}}{j!(3)_j(n^2+n-2-j)!}\;\right.\\ 
    \Biggl.\times _3F_2\left(n,n+2,n^2+n+1;n+1,n^2+n-1-j;\frac{c_\eta(z-n)}{z-c_\eta}\right) \Biggr) H(z-n).
  \end{multline}
\end{corollary}

The expression corresponding to $\eta=0$ agrees with a previous result given in \cite[Theorem 3.3]{ref:prathapSIAM} as shown in the following corollary.
\begin{corollary}
  \label{corretazero}
  The exact p.d.f. of $\kappa_{\text{SC}}^2(\mathbf{X})$ given in Theorem \ref{thm exact pdf} simplifies, for $\eta=0$, giving
  \begin{align}
  \label{eq dem pdf white}
      f^\alpha_{\kappa_{\text{SC}}^2(\mathbf{X})}(z)&= \Gamma\left(mn\right)\left(\prod_{j=0}^\alpha\frac{n+j}{(j+1)!}\right)(z-n)^{mn-\alpha-2}z^{-mn} \nonumber \\
  &\qquad\qquad\qquad\qquad \times \sum_{k_1=0}^{n+\alpha-2} \text{...} \sum_{k_{\alpha}=0}^{n-1}\left(\prod_{j=1}^{\alpha}(-1)^{k_j}\frac{(-n-\alpha+j+1)_{k_j}}{(j+1)_{k_j}k_j!} (z-n)^{-k_j}\right) \nonumber \\
  &\qquad\qquad\qquad\qquad\qquad\qquad\qquad\qquad\quad \times \frac{\Delta_\alpha(\textbf{c})}{\Gamma\left(mn-\alpha-1-\sum_{j=1}^{\alpha}k_j\right)}H(z-n) \text{ }
  \end{align}
  which coincides with \cite[Theorem 3.3]{ref:prathapSIAM}. Here, $\textbf{c}=\{c_1,c_2,\dots,c_\alpha\}$ with $c_j=j+k_j$ and $m=n+\alpha$.
\end{corollary}
\begin{IEEEproof}
  See Appendix \ref{appendix_corretazero}.
\end{IEEEproof}

Interestingly, capitalizing on the joint eigenvalue density (\ref{cor joint}), we can derive novel expressions for the p.d.f. and the c.d.f. of the minimum eigenvalue $\lambda_1$ of the single-spiked Wishart-Laguerre ensemble as shown in the following lemma.
\begin{lemma}
  The p.d.f. and c.d.f. of the minimum eigenvalue $\lambda_1$ of the single-spiked Wishart-Laguerre ensemble is given, respectively, by
  \begin{equation}
    f_{\lambda_1}^\alpha(x) = \frac{(n-1)! x^\alpha e^{-x\left(n-c_\eta\right)}}{(n+\alpha-1)!(\eta+1)^{\alpha}} \det\left[\left(n+i-c_\eta\right)\left(-\eta\right)^i \;\; L_{n+i-j-1}^{(j+1)}(-x) \right]_{\substack{i=0,...,\alpha \\ j=1,...,\alpha}},
  \end{equation}
  \begin{equation}
    F_{\lambda_1}^\alpha(x) = 1- \frac{e^{-x\left(n-c_\eta\right)}}{(\eta+1)^{\alpha}} \det \left[(-\eta)^i \;\;\; L_{n+i-j}^{(j-1)}(-x)\right]_{\substack{i=0,...,\alpha \\ j=1,...,\alpha}} \label{cdf lammin}.
  \end{equation}
\end{lemma}
\begin{IEEEproof}
  See Appendices \ref{appendix_pdfmineig} and \ref{appendix_cdfmineig}.
\end{IEEEproof}

Similar results (i.e., $\alpha$ dependent determinant size) appear in the literature \cite{ref:prathapSIAM, ref:Forrester, ref:akemann}  for $\lambda_1$ of Wishart-Laguerre ensemble (i.e., without the single-spiked covariance or $\boldsymbol{\Sigma}=\mathbf{I}_n$). For a fully correlated Wishart matrix, a more complicated $\alpha$ dependent result for the minimum eigenvalue has been reported in \cite{ref:tim}. However, for single-spiked case, the above results are much compact and easy to handle. As a sanity check, for $\alpha=0$, we obtain $F_{\lambda_1}^0 (x)=1-e^{-x\left(n-c_\eta\right)}$, which coincides with \cite[Eq. 2.15]{ref:petermin}, \cite[Eq. 3.23]{ref:peterBeta}, and \cite[Eq. 21]{ref:rathnasiam}.

\begin{figure}[tbhp]
    \centering
    \includegraphics[width=12cm]{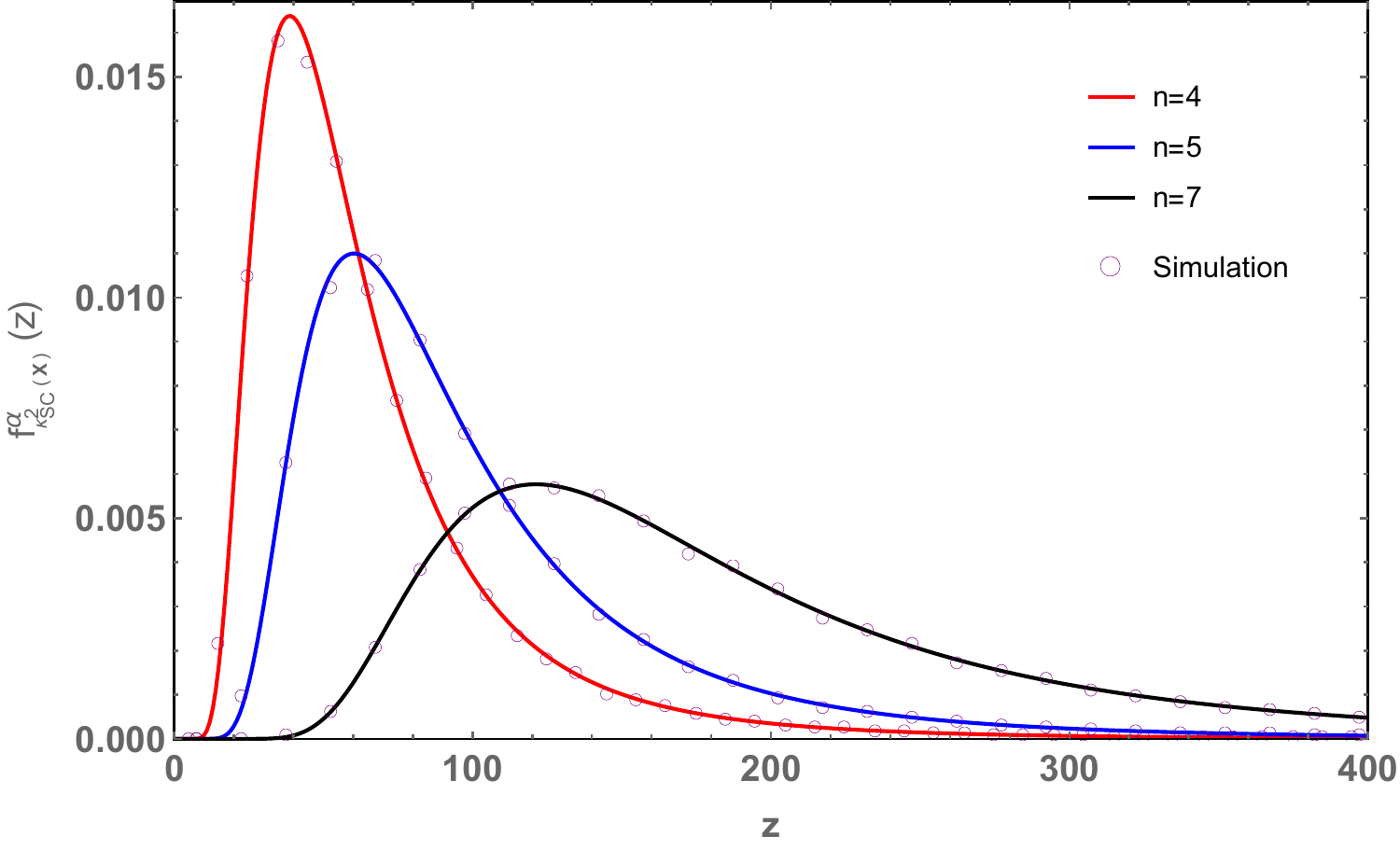}
    \caption{Comparison of simulated data points and the analytical p.d.f. $f^\alpha_{\kappa^2_{\rm SC}(\mathbf{X})}(z)$ for different values of $n$ with $\alpha=3$ and $\eta=10$.}
    \label{fig_different_n}
\end{figure}

\begin{figure}[tbhp]
    \centering
    \includegraphics[width=12cm]{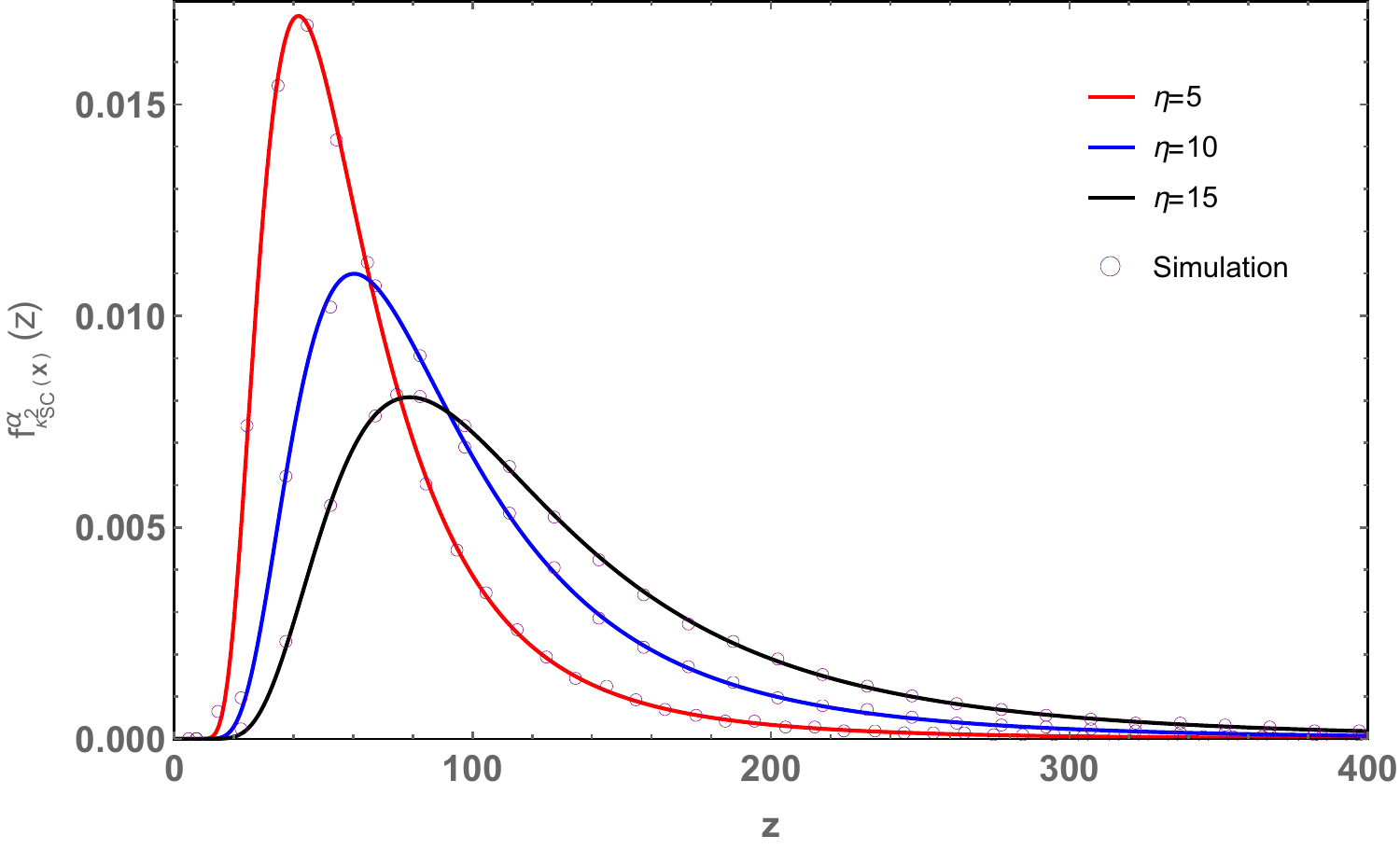}
    \caption{Comparison of simulated data points and the analytical p.d.f. $f^\alpha_{\kappa^2_{\rm SC}(\mathbf{X})}(z)$ for different values of $\eta$ with $n=5$ and $\alpha=3$.}
    \label{fig_different_eta}
\end{figure}

\begin{figure}[tbhp]
    \centering
    \includegraphics[width=12cm]{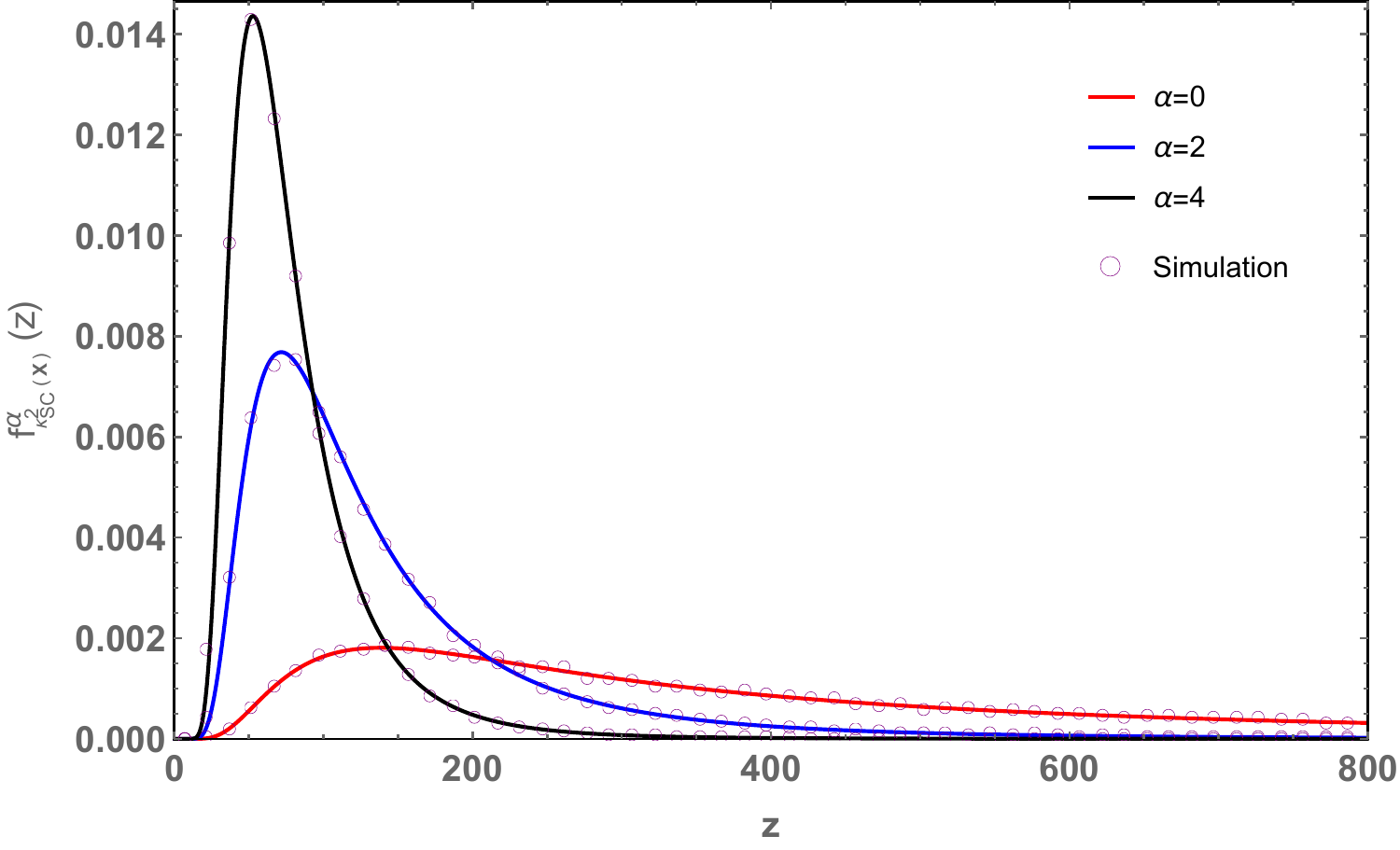}
    \caption{Comparison of simulated data points and the analytical p.d.f. $f^\alpha_{\kappa^2_{\rm SC}(\mathbf{X})}(z)$ for different values of $\alpha$ with $n=5$ and $\eta=10$.}
    \label{fig_different_alpha}
\end{figure}

Figures \ref{fig_different_n} and \ref{fig_different_eta} compare the analytical p.d.f. result of $\kappa_{\text{SC}}^2(\mathbf{X})$ for the deformed Laguerre-Wishart matrix computed based on Theorem \ref{thm exact pdf} with simulated data. In particular, Fig. \ref{fig_different_n} depicts the p.d.f. of $\kappa_{\text{SC}}^2(\mathbf{X})$ for various $n$ and $m$ configurations with $\eta=10$. Figure \ref{fig_different_eta} shows the effect of $\eta$ on the p.d.f. for the fixed configuration  $n=5$ and $m=8$ (i.e., $\alpha=3$). Moreover, Fig. \ref{fig_different_alpha} shows the effect of $m$ on the p.d.f. of $\kappa_{\text{SC}}^2(\mathbf{X})$ for fixed $\eta$ and $n$. As depicted in the figure, as $m$ increases for fixed $n$ and $\eta$ (i.e., $n=5$ and $\eta=10$), the p.d.f. of $\kappa_{\text{SC}}^2(\mathbf{X})$ tends to concentrate towards $z=16$. To further examine this behaviour, let us focus on the limiting behaviour of $\kappa_{\text{SC}}^2(\mathbf{X})$ as $m\to\infty$. Noting that $\lambda_1/m$ tends to $1$ almost surely and $\text{tr}(\mathbf{W})/m=\sum_{j=1}^n \lambda_j/m$ tends to $\eta+n+1$ in probability\footnote{It can easily be proved that the characteristic function of $\text{tr}(\mathbf{W})/m$ given by $\mathbb{E}\left\{e^{-j\omega\text{tr}(\mathbf{W})/m }\right\}=\frac{1}{\left(1+(1+\eta)j\omega/m\right)^m (1+j\omega/m)^{m(n-1)}}$ converges to $e^{-j\omega(n+\eta+1)}$ as $m\to\infty$ for fixed $n$. Therefore, following the L\'evy's continuity theorem, we obtain $\text{tr}(\mathbf{W})/m$ converges to $n+\eta+1$ weakly (i.e., converges in distribution). The final claim follows by noting that weak convergence to a constant amounts to convergence in probability to the same constant.}, following continuous mapping theorem \cite{ref:vaart}, we can conclude that $\kappa_{\text{SC}}^2(\mathbf{X})=\left(\sum_{j=1}^n \lambda_j/m\right)/\left(\lambda_1/m\right)$ tends to $\eta+n+1$ in probability. This explains reason behind leftward shift of the p.d.f.s in Fig. \ref{fig_different_alpha} as $\alpha$ (i.e., $m$) increases for fixed $\eta$.
 
{\color{blue}Having statistically characterized the p.d.f. of $\kappa_{\text{SC}}^2(\mathbf{X})$, we now focus on developing the ROC curves for the test statistic $T(\boldsymbol{\lambda})$ in (\ref{Teststat}). Since the dependency of the p.d.f. of $\kappa_{\text{SC}}^2(\mathbf{X})$ on the perturbation power $\eta$ is of  paramount importance in the sequel, we rewrite $f^\alpha_{\kappa_{\text{SC}}^2(\mathbf{X})}(z)$ as $f^\alpha_{\kappa_{\text{SC}}^2(\mathbf{X})}(z,\eta)$ to indicate it. Therefore, under this setting, the detection\footnote{This is also known as the power of the test.} and false alarm probabilities can be written with the help of (\ref{Teststat}), (\ref{CRsetting}), and Theorem \ref{thm exact pdf} as
\begin{align}
    P_D^{\alpha}(\eta_{\rm snr}, \xi_{\rm th})&=\Pr\left\{T(\boldsymbol{\lambda})>\xi_{\rm th}|\mathcal{H}_1\right\},\label{pd}\\
    P_F^{\alpha}(\xi_{\rm th})&= \Pr\left\{T(\boldsymbol{\lambda})>\xi_{\rm th}|\mathcal{H}_0\right\}\label{pf}
\end{align}
where $\xi_{\rm th}$ denotes a certain threshold\footnote{Since $T(\boldsymbol{\lambda})$ is a continuous random variable, the threshold $\xi_{\rm th}$ is chosen such that $P_F^{\alpha}(\xi_{\rm th})=\delta\in(0,1)$.}  and 
\begin{align}
\label{Tdistri}
    T(\boldsymbol{\lambda})\sim\left\{
    \begin{array}{ll}
    f^\alpha_{\kappa_{\text{SC}}^2(\mathbf{X})}(z,\eta_{\rm snr}) & \text{under $\mathcal{H}_1$}\\
    f^\alpha_{\kappa_{\text{SC}}^2(\mathbf{X})}(z,0) & \text{under $\mathcal{H}_0$}
    \end{array}\right.
\end{align}
with $\eta_{\rm snr}=\displaystyle\frac{\gamma ||\mathbf{h}||^2}{\sigma^2}$ denoting the SNR. Now the $(P^\alpha_D, P^\alpha_F)\in [0,1]^2 $ pair characterizes   the detector  and is called as the  ROC profile. In general, obtaining an explicit functional relationship between $P^\alpha_D$ and $P^\alpha_F$ by eliminating the parametric dependency on $\xi_{\rm th}$ seems an arduous task. Nevertheless, when $\alpha$ admits zero, an explicit relationship between them is possible. To show this, let us first  
write $P_F^\alpha$ in (\ref{pf}) as
\begin{align}
    P_F^{\alpha}(\xi_{\rm th})=\int_{\xi_{\rm th}}^\infty 
    f^\alpha_{\kappa_{\text{SC}}^2(\mathbf{X})}(z,0) {\rm d}z
\end{align}
which can be further simplified in view of (\ref{eq dem pdf white}) with $\alpha=0$ to obtain
\begin{align}
    P_F^{0}(\xi_{\rm th})=n(n^2-1)
    \int_{\xi_{\rm th}}^\infty
    \frac{(z-n)^{n^2-2}}{z^{n^2}} {\rm d}z=1-\left(1-\frac{n}{\xi_{\rm th}}\right)^{n^2-1}.
\end{align}
Therefore, we get \cite[Eq. 15]{ref:ayse}
\begin{align}
\label{thresh}
    \xi_{\rm th}=\frac{n}{1-\left(1-P^0_F\right)^{\frac{1}{n^2-1}}}.
\end{align}
Now following (\ref{pd}), we write the detection probability as
\begin{align}
    P_D^{\alpha}(\eta_{\rm snr}, \xi_{\rm th})=\int_{\xi_{\rm th}}^\infty
    f^\alpha_{\kappa_{\text{SC}}^2(\mathbf{X})}(z,\eta_{\rm snr}) {\rm d}z
\end{align}
from which, in view of (\ref{eq_pdf_alpha0}), we obtain
\begin{align}
\label{powerint}
   P_D^{0}(\eta_{\rm snr}, \xi_{\rm th})=
   \frac{n(n^2-1)}{(\eta_{\rm snr}+1)^n}
   \int_{\xi_{\rm th}}^\infty  \frac{(z-n)^{n^2-2}}{ (z-c_{\eta_{\rm snr}})^{n^2}}  {}_3F_2\left(n-1,n+1,n^2;n,n^2-1;\frac{c_{\eta_{\rm snr}}(z-n)}{z-c_{\eta_{\rm snr}}}\right) {\rm d}z.
\end{align}
To facilitate further analysis, we use (\ref{3F2decom}) with $\alpha=0$ to rewrite
${}_3F_2$ as
\begin{align}
    {}_3F_2\left(n-1,n+1,n^2;n,n^2-1;\frac{c_{\eta_{\rm snr}}(z-n)}{z-c_{\eta_{\rm snr}}}\right)
    &=\frac{(z-c_{\eta_{\rm snr}})^{n+1}}{(n+1)\left[(1-c_{\eta_{\rm snr}})z+(n-1)c_{\eta_{\rm snr}}\right]^{n+1}}\nonumber\\
    &\qquad + \frac{z (n-c_{\eta_{\rm snr}}) (z-c_{\eta_{\rm snr}})^{n-1}}{(n+1)\left[(1-c_{\eta_{\rm snr}})z+(n-1)c_{\eta_{\rm snr}}\right]^{n}},
\end{align}
which upon substituting into (\ref{powerint}) followed by some algebraic manipulation gives
\begin{align}
\label{powerintmid}
    P_D^{0}(\eta_{\rm snr}, \xi_{\rm th})&=
    n(n-1)(\eta_{\rm snr}+1) \int_{\xi_{\rm th}}^\infty \frac{(z-n)^{n^2-2}}{(z-c_{\eta_{\rm snr}})^{n^2-n-1} \left[z+(n-1)\eta_{\rm snr}\right]^{n+1}} \;{\rm d}z \nonumber\\
    &\qquad +n(n-1)(n-c_{\eta_{\rm snr}}) \int_{\xi_{\rm th}}^\infty \frac{z(z-n)^{n^2-2}}{(z-c_{\eta_{\rm snr}})^{n^2-n+1} \left[z+(n-1)\eta_{\rm snr}\right]^{n}}\; {\rm d}z.
\end{align}
{\color{blue}
Now we may apply the substitution $x=z+(n-1)\eta_{\rm snr}$ into (\ref{powerintmid}) and evaluate the resultant integrals to obtain the probability of detection as
\begin{align}
\label{powerfinal}
    P_D^{0}(\eta_{\rm snr}, \xi_{\rm th})&=
    \frac{n(n-1)\left(1+\eta_{\rm snr}\right)^{n^2}}{\eta_{\rm snr}^{n^2-1} \left(n+(n-1)\eta_{\rm snr}\right)}
    \sum_{k=0}^{n^2-2}
    \sum_{\ell=0}^2(-1)^k
    \frac{ \binom{n^2-2}{k}(n-1)^{\varepsilon_\ell}}{(\eta_{\rm snr}+1)^{k+\ell}}
    \mathcal{P}_{k+\ell}(\xi_{\text{th}},\eta_{\rm snr})
\end{align}
where 
\begin{align}
    \mathcal{P}_{M}(z,\eta_{\rm snr})= \left\{\begin{array}{ll}
    \displaystyle \frac{1}{n-M}\left[1- 
    \left(1-\frac{\eta_{\rm snr} [n+(n-1)\eta_{\rm snr}]}{[\eta_{\rm snr}+1][z+(n-1)\eta_{\rm snr}]}\right)^{n-M}\right] & \text{for $n\neq M$}\\
    \displaystyle -\ln\left(1-\frac{\eta_{\rm snr} [n+(n-1)\eta_{\rm snr}]}{[\eta_{\rm snr}+1][z+(n-1)\eta_{\rm snr}]}\right) & \text{for $n=M$}
    \end{array}\right.,
\end{align}
$\binom{n}{k}=n!/k!(n-k)!$ is the binomial coefficient, and $\varepsilon_\ell=\frac{1}{2}(1-(-1)^\ell)$. Consequently, (\ref{thresh}) and (\ref{powerfinal}) together determine the ROC curve of the test statistic $T$ corresponding to $\alpha=0$ for an arbitrary SNR $\eta_{\rm snr}$.

Figure \ref{fig_roc_alpha0_dietectionff_eta} depicts the analytical and simulated ROC curves for different SNR values with $\alpha=0$ and $n=4$. The positive effect of SNR on the probability of detection is clearly visible from the figure. To further highlight this and to investigate the effect of the sample size (i.e., $m$) on the probability of detection, in Fig. \ref{fig_roc_diff_alpha}, we plot the detection probability versus false alarm probability for various values of $m$ and $\eta_{\rm snr}$ with $n=5$. Here the theoretical ROC curves corresponding to $\alpha=0$ have been generated by using (\ref{thresh}) and (\ref{powerfinal}), whereas, for $\alpha\neq 0$, we have numerically integrated the analytical p.d.f.s given in  (\ref{eq_pdf_Z}) and (\ref{eq dem pdf white}). Again, as expected, both parameters $m$ and $\eta_{\rm snr}$ affect the detection probability positively. Having understood the joint effect of the latter parameters, let us now focus on demonstrating the joint effect of the number of samples and the number of receive antennas (i.e., $n$) or equivalently the number of single antenna secondary users on the detection probability. To this end, Fig. \ref{fig_roc_diff_n} depicts the  ROC curves corresponding to different values of $\alpha$ for various configurations of $m$ and $n$ such that $\alpha=m-n$ is fixed with $\mathbf{h}=(1 \; 1\; \ldots\;1)^T$. It is noteworthy that this particular choice of $\mathbf{h}$ gives $||\mathbf{h}||^2=n$, whereas, for $\mathbf{h}\sim\mathcal{CN}_n(\mathbf{0},\mathbf{I}_n)$ (i.e., Rayleigh fading channels), we obtain the almost sure limit $\displaystyle \lim_{n\to\infty }\frac{||\mathbf{h}||^2}{n}\to 1$. This in turn suggests that, for large enough $n$, $||\mathbf{h}||^2$ can be approximated with $n$; thereby drawing an analogy between the two channels. As can be seen from the figure, 
 increasing both $m$ and $n$ with their difference fixed leads to an improved detection probability for an arbitrary transmit SNR (i.e., $\gamma/\sigma^2$). The reason behind this improvement for fixed $\alpha$ is the  dependency of $\eta_{\rm snr}$ on $n$. To be specific, for each $n$ with fixed $\alpha$, we have $\eta_{\rm snr}=\gamma n/\sigma^2$; therefore, $\eta_{\rm snr}$ grows linearly with $n$. In the setting of Fig. \ref{fig_roc_diff_n}, the values of $\eta_{\rm snr}$ corresponding to $n=2,5,8$ are given, respectively, by $\eta_{\rm snr}\approx 8, 12, 14$ dB. }

\begin{figure}[tbhp]
    \centering
    \includegraphics[width=12cm]{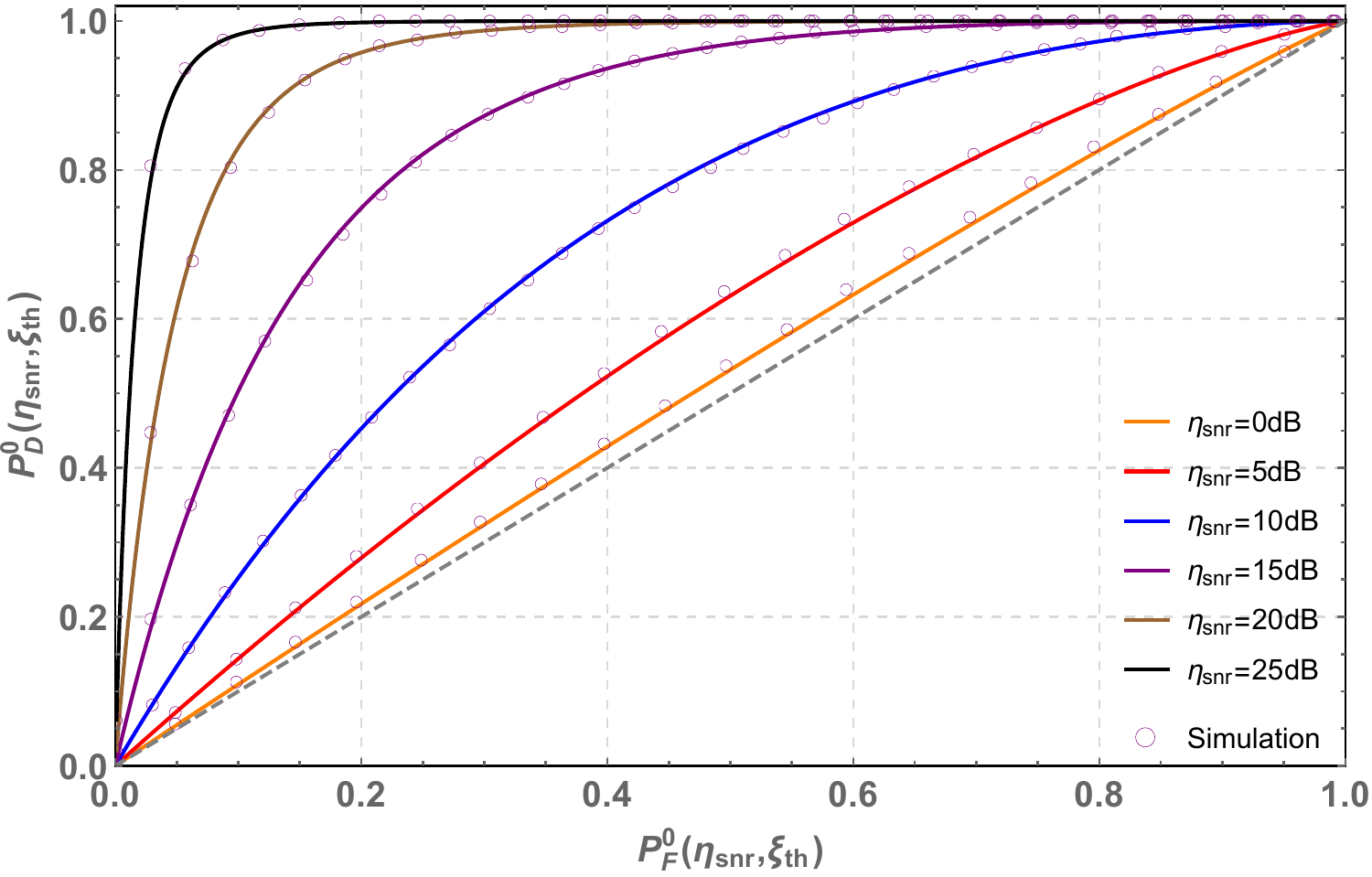}
    \caption{Probability of detection versus false alarm probability; simulation and analytical result comparison. ROC curves are shown for different $\eta_{\rm snr}$ with $\alpha=0$ and $n=4$}
    \label{fig_roc_alpha0_dietectionff_eta}
\end{figure}

\begin{figure}[tbhp]
    \centering
    \includegraphics[width=12cm]{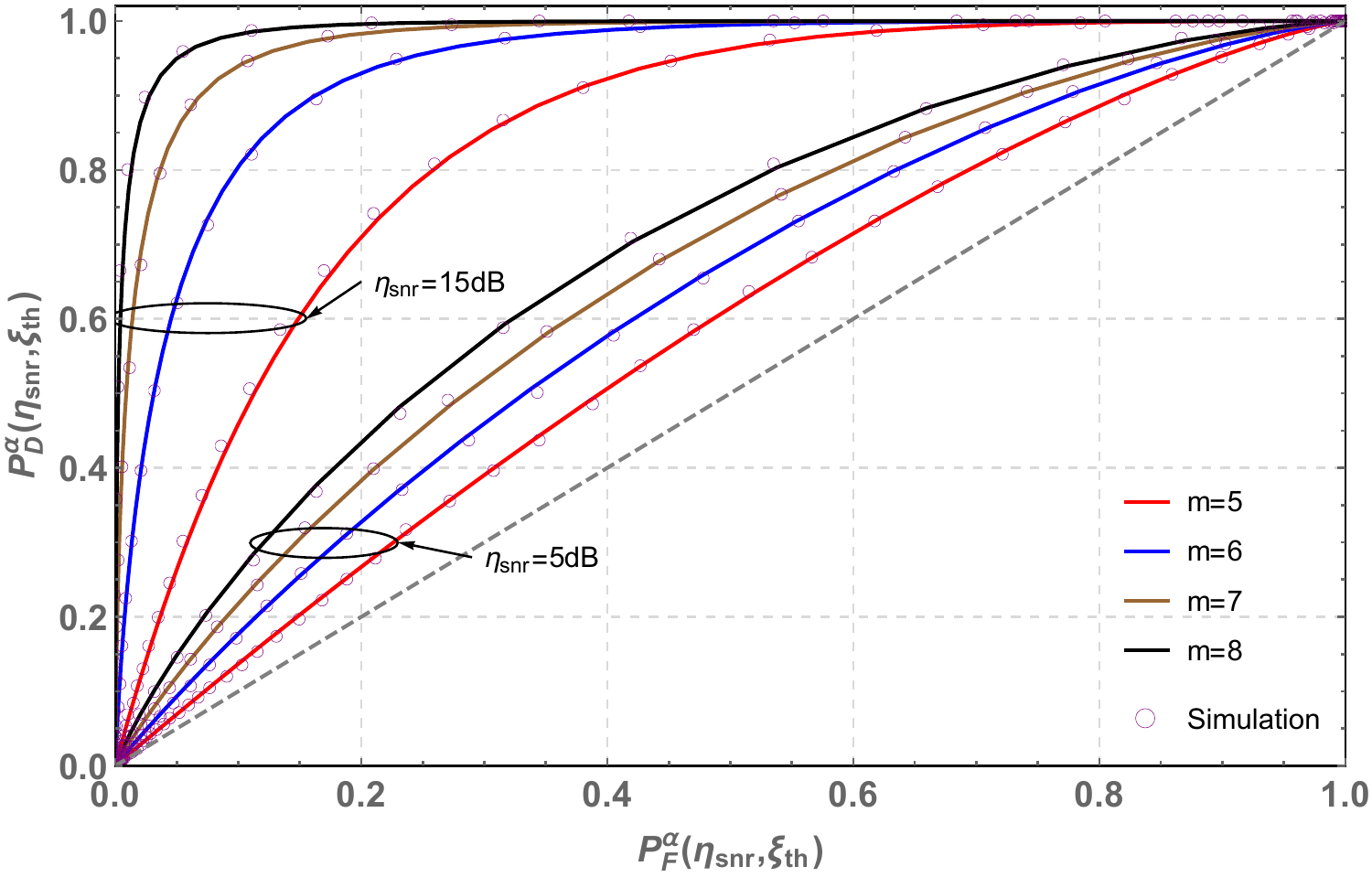}
    \caption{Probability of detection versus false alarm probability; simulation and analytical result comparison. ROC curves are shown for different values of $m$ with  $n=5$ and $\eta_{\rm snr}=5,15\; \text{dB}$.}
    \label{fig_roc_diff_alpha}
\end{figure}



\begin{figure}[tbhp]
    \centering
    \includegraphics[width=12cm]{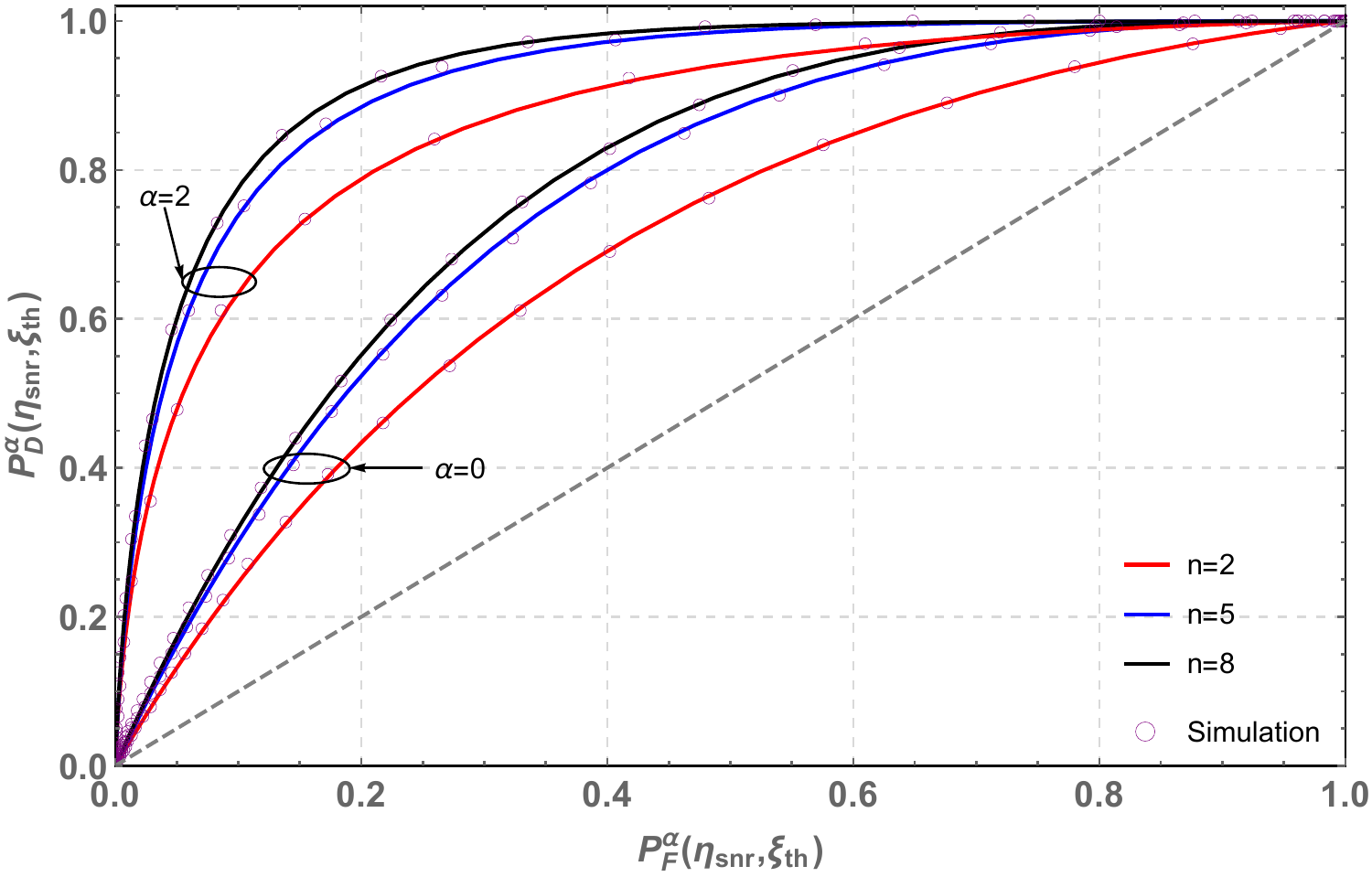}
    \caption{Probability of detection versus false alarm probability; simulation and analytical result comparison. ROC curves are shown for different $n$ and $\alpha$ with the transmit  ${\rm SNR}=5$ dB.}
    \label{fig_roc_diff_n}
\end{figure}

}

\section{Asymptotic Characterization of $\kappa_{\text{SC}}^2(\mathbf{X})$}

{\color{blue}In this section, we investigate the distribution of $\kappa_{\text{SC}}^2(\mathbf{X})$ over two pertinent asymptotic regimes. In particular, our focus is on the two regimes: $m,n\to\infty$ such that $m-n$ is fixed and $m,n\to\infty$ such that $n/m\to c \in (0,1)$. Our general strategy is to first characterize the limiting distributions (i.e., weak limits) of suitably centered and scaled $\lambda_1$ and to obtain the weak limits of suitably centered and scaled $\kappa_{\text{SC}}^2(\mathbf{X})$ corresponding to the above two regimes by invoking the Slutsky's lemma \cite{ref:vaart} subsequently.} 

\subsection{Limiting Distribution in the Fixed $m-n$ Regime}
In this subsection, we use the exact c.d.f. expression (\ref{cdf lammin}) to investigate the c.d.f. of $\kappa_{\text{SC}}^2(\mathbf{X})$, suitably scaled, for fixed $\alpha$ when $m,n\to\infty$.  As such, we have the following theorem.
\begin{theorem}\label{thm asy lam}
    The scaled random variable $X_n=\mu n \lambda_{1}$ with $\eta$ scaled as $\eta=\frac{\rho}{n}$, where $\mu, \rho \in \mathbb{R}^+$ are arbitrary constants, converges in distribution to another random variable $X$ with the following c.d.f as $m,n \rightarrow \infty$ with $\alpha=m-n$ fixed:
    \begin{align}
        F_X^\alpha(x) &= 1-e^{-\frac{x}{\mu}}\det\left[I_{j-i}\left(2\sqrt{\frac{x}{\mu}}\right)\right]_{i,j=1,\dotsc,\alpha}
    \end{align}
    where $I_k(\cdot)$ denotes the modified Bessel function of the second kind and order $k$.
\end{theorem}
\begin{IEEEproof}
We may use (\ref{cdf lammin}) and the the definition of Laguerre polynomial given in (\ref{lagdef}) to obtain
\begin{align}
    \text{Pr}\left(\lambda_{1}\geq x\right) &= \frac{e^{-x\left(n-c_\eta\right)}}{(\eta+1)^\alpha} \det\left[(-\eta)^i \;\;\; \frac{(j)_{n+i-j}}{(n+i-j)!}\sum_{k_j=0}^{n+i-j}\frac{(-n-i+j)_{k_j}}{(j)_{k_j}}\frac{(-x)^{k_j}}{k_j!}\right]_{\substack{i=0,\dotsc,\alpha \\ j=1,\dotsc,\alpha}}.
\end{align} 
Further manipulation of the above determinant seems an arduous task due to the $i,j$-dependent summation upper limits. To circumvent this difficulty, 
noting that $(-n-i+j)_{k_j} = 0$ for $k_j>n+i-j$, the above expression can be re-written as
\begin{align}
\label{eq cdf mani}
    \text{Pr}\left(\lambda_{1}\geq x\right) &= \frac{e^{-x\left(n-c_\eta\right)}}{(\eta+1)^\alpha} \det\left[(-\eta)^i \;\;\; \frac{(n+i-1)!}{(j-1)!(n+i-j)!}\sum_{k_j=0}^{n+\alpha-j}\frac{(-n-i+j)_{k_j}}{(j)_{k_j}}\frac{(-x)^{k_j}}{k_j!}\right]_{\substack{i=0,\dotsc,\alpha \\ j=1,\dotsc,\alpha}} \text{ .}
\end{align}
Again, to eliminate the $i,j$-dependency in the numerator of each summation, we use the decomposition,
\begin{align}
    (-n-i+j)_{k_j} &= (-n-i+j)_{k_j}\frac{(-n-\alpha+j)_{k_j}}{(-n-\alpha+j)_{k_j}} \\
    &= \frac{(n+i-j)!}{(n+\alpha-j)!}(-n-\alpha+j)_{k_j}\prod_{p=0}^{\alpha-i-1}(\hat{c}_j-p)
\end{align}
where $\hat{c}_j = n+\alpha-j-k_j$, in (\ref{eq cdf mani})  with some algebraic manipulations to obtain
\begin{align}
    \text{Pr}\left(\lambda_{1} \geq x\right) &= \frac{e^{-x\left(n-c_\eta\right)}}{(\eta+1)^\alpha} \prod_{j=1}^{\alpha}\frac{1}{(j-1)!} \sum_{k_1=0}^{n+\alpha-1}\dotsc\sum_{k_\alpha=0}^n \left(\prod_{j=1}^\alpha \frac{(-n-\alpha+j)_{k_j}}{(j)_{k_j}}\frac{(-x)^{k_j}}{k_j!}\right) \nonumber \\
    &\qquad\qquad\qquad\qquad\qquad \times \det\left[\frac{(-\eta)^i(n-1)!}{(n+i-1)!}\;\;\; \prod_{p=0}^{\alpha-i-1}(\hat{c}_j-p)\right]_{\substack{i=0,\dotsc,\alpha \\ j=1,\dotsc,\alpha}}
\end{align}
where an empty product is interpreted as unity.
Now, to determine the limiting distribution, we consider the scaled random variable $\mu n\lambda_1$ and the scaled parameter $\eta=\rho/n$. As such, after using some elementary  limiting arguments we arrive at
\begin{align}
\label{eq lim cdf}
    \lim_{n\rightarrow\infty}\text{Pr}\left(\mu n \lambda_{1} \geq x\right) &= \frac{e^{-\frac{x}{\mu}}}{\prod_{j=1}^\alpha(j-1)!}\sum_{k_1=0}^\infty\dotsc\sum_{k_\alpha=0}^\infty\left(\prod_{j=1}^\alpha\frac{1}{(j)_{k_j}k_j!}\frac{x^{k_j}}{\mu^{k_j}}\right)  \lim_{n\rightarrow\infty}
    \Omega(\rho, \alpha,n)
\end{align}
where 
\begin{align*}
    \Omega(\rho, \alpha,n)=
    \det\left[\frac{(-\rho)^i(n-1)!}{n^i(n+i-1)!} \;\;\; \prod_{p=0}^{\alpha-i-1}(\hat{c}_j-p)\right]_{\substack{i=0,\dotsc,\alpha \\ j=1,\dotsc,\alpha}}.
\end{align*}
To facilitate further analysis, we need to obtain the limiting value of the above determinant. To this end, we use some algebraic manipulations to yield
\begin{align}
    \Omega(\rho, \alpha, n) &= \det\left[ \frac{(-1)^i\rho^i}{n^i\prod_{q=0}^{i-1}(n+q)} \;\;\; \prod_{p=0}^{\alpha-i-1}(\hat{c}_j-p)\right]_{\substack{i=0,\dotsc,\alpha \\ j=1,\dotsc,\alpha}}
\end{align}
from which we obtain 
upon invoking \cite[Lemma A.1]{ref:prathapSIAM}
\begin{align}
    \Omega(\rho, \alpha, n) &= \left[
        \begin{matrix}
            (-1)^0 + o\left(\frac{1}{n^2}\right) & \hat{c}_1^\alpha & \dotsc & \hat{c}_\alpha^\alpha \\
            \frac{(-1)^1 \rho}{n^2}+o\left(\frac{1}{n^4}\right) & \hat{c}_1^{\alpha-1} & \dotsc & \hat{c}_\alpha^{\alpha-1} \\
            \vdots & \vdots & & \vdots \\
            \frac{(-1)^i \rho^i}{n^i\prod_{q=0}^{i-1}(n+q)}+o\left(\frac{1}{n^{2(i+1)}}\right) & \hat{c}_1^{\alpha-i} & \dotsc & \hat{c}_\alpha^{\alpha-i} \\
            \vdots & \vdots & & \vdots \\
            \frac{(-1)^\alpha\rho^\alpha}{n^\alpha\prod_{q=0}^{\alpha-1}(n+q)} & 1 & \dotsc & 1
        \end{matrix}
    \right]
\end{align}
where $o(\cdot)$ denotes the little-$o$ notation.
Now following \cite{ref:prathapJMVA}, we obtain
\begin{align*}
    \lim_{n\to\infty} \Omega(\rho, \alpha,n)=\Delta_\alpha\left({\mathbf{c}}\right)
\end{align*}
where ${\mathbf{c}}=({c}_1,{c}_2,\ldots,{c}_{\alpha})$with $c_\ell=\ell+k_\ell$.
Consequently, (\ref{eq lim cdf}) assumes the form
\begin{align*}
    \lim_{n\rightarrow\infty}\text{Pr}\left(\mu n \lambda_{1} \geq x\right) &= \frac{e^{-\frac{x}{\mu}}}{\prod_{j=1}^\alpha(j-1)!}\sum_{k_1=0}^\infty\dotsc\sum_{k_\alpha=0}^\infty\left(\prod_{j=1}^\alpha\frac{1}{(j)_{k_j}k_j!}\frac{x^{k_j}}{\mu^{k_j}}\right) \Delta_\alpha\left({\mathbf{c}}\right)
\end{align*}
from which we obtain using \cite{ref:prathapJMVA} 
\begin{align}
    \lim_{n\rightarrow\infty}\text{Pr}\left(\mu n\lambda_{1} \geq x\right) &= e^{-\frac{x}{\mu}}\det\left[I_{j-i}\left(2\sqrt{\frac{x}{\mu}}\right)\right]_{i,j=1,\dotsc,\alpha}.
\end{align}
Finally, we make use of the relation
\begin{align}
    F_{X}^\alpha(x) &= \lim_{n\rightarrow\infty}F_{\mu n \lambda_{1}}^\alpha(x) = 1-\lim_{ n\rightarrow\infty}\text{Pr}\left(\mu n\lambda_{1} \geq x\right)
\end{align}
to conclude the proof.
\end{IEEEproof}

Having asymptotically characterized $\lambda_1$, we are now in a position to present the weak limit of properly scaled $\kappa_{\text{SC}}^2(\mathbf{X})$ which is given by the following corollary.

\begin{corollary}\label{cor asy cdf}
The scaled random variable $V_n=\kappa_{\text{SC}}^2(\mathbf{X})/\mu n^3$ with $\eta=\frac{\rho}{n}$, where $\mu, \rho \in \mathbb{R}^+$ are arbitrary constants, converges in distribution to another random variable $V$ with the following c.d.f. and p.d.f., respectively, as $m,n\xrightarrow{}\infty$ with $\alpha=m-n$ fixed:
\begin{align}
    F_V^\alpha(v)&=e^{-\frac{1}{\mu v}} \det \left[I_{j-i}\left(\frac{2}{\sqrt{\mu v}}\right)\right]_{i,j=1,..,\alpha}H(v)\label{eq dem asym}\\
    f_V^\alpha(v)&=\frac{e^{-\frac{1}{\mu v}}}{\mu v^2} \det \left[I_{j-i+2}\left(\frac{2}{\sqrt{\mu v}}\right)\right]_{i,j=1,..,\alpha} H(v)\label{asy pdf theo}.
\end{align}
\end{corollary}
\begin{IEEEproof}
Since, from Theorem \ref{thm asy lam}, as $m,n\to\infty$ with $\alpha$ fixed, $\mu n \lambda_1$ converges in distribution to $X$ and $\sum_{j=1}^n \lambda_j/n^2$ converges in probability to $1$\footnote{It can be proved that the characteristic function of $\sum_{j=1}^n\lambda_j/n^2=\text{tr}(\mathbf{W})/n^2$ given by $\mathbb{E}\left\{e^{-j\omega\text{tr}(\mathbf{W})/n^2 }\right\}=\frac{1}{\left(1+(1+\rho/n)j\omega/n^2\right)^{(n+\alpha)} (1+j\omega/n^2)^{(n+\alpha)(n-1)}}$ converges to $e^{-j\omega}$ as $n\to\infty$ for fixed $\alpha$. Therefore, following the L\'evy's continuity theorem, we obtain $\text{tr}(\mathbf{W})/n^2$ converges to $1$ weakly (i.e., converges in distribution). The final claim follows by noting that weak convergence to a constant amounts to convergence in probability to the same constant.}, we may use the Slutsky's lemma \cite{ref:vaart} to obtain the weak limit of $\mu n^3/\kappa_{\text{SC}}^2(\mathbf{X})$. Consequently, we invoke the continuous mapping theorem \cite{ref:vaart} to conclude the proof of (\ref{eq dem asym}).

{\color{blue}Whereas the derivative of $F_V^\alpha(v)$ with respect to $v$ gives, in principle, the p.d.f. of $V$,  that approach does not yield any simple expression for the p.d.f. To overcome this difficulty, here we directly evaluate the limiting p.d.f. of $V$ starting from the p.d.f. of $\kappa_{\text{SC}}^2(\mathbf{X})$ given in Theorem (\ref{thm exact pdf}). To this end, we multiply and divide (\ref{eq_pdf_Z}) by the factor $(n+\alpha-j-k_j-1)!$ with some algebraic manipulation to yield
\begin{align}
    f_{\kappa_{\text{SC}}^2(\mathbf{X})}^\alpha(z)&=\frac{H(z-n)}{(\eta+1)^{n+\alpha}} \sum_{k_1=0}^{n+\alpha-2}\text{...}\sum_{k_\alpha=0}^{n-1} \left(\prod_{j=1}^\alpha \frac{(n+\alpha-j-1)!}{(n+\alpha-j-k_j-1)!(j+k_j+1)!k_j!} \right) \nonumber \\
     & \times \frac{(z-n)^{(n-1)(n+\alpha+1)-\sum_{j=1}^\alpha k_j-1}}{(z-c_\eta)^{n(n+\alpha)}} \det \left[(n+\alpha)!\mathcal{G}_i(z,\eta) \;\; \prod_{\ell=0}^{\alpha-i-1}(\Tilde{c_j}-\ell)\right]_{\substack{i=0,...,\alpha \\ j=1,...,\alpha} }
\end{align}
where $\tilde{c_j}=n+\alpha-j-k_j-1$. Now noting the fact that
 $f_V^\alpha(v)=\displaystyle \lim_{n\xrightarrow{}\infty}\mu n^3 f_{\kappa_{\text{SC}}^2(\mathbf{X})}^\alpha(\mu n^3 v)$ and $\eta=\rho/n$, we apply the fundamental limiting arguments as $n\to\infty$ followed by some algebraic manipulation to arrive at
 \begin{align}
 \label{eq_asymp_substitution1}
     f_V^\alpha(v) &= \lim_{n\xrightarrow{}\infty}\mu n^3 f_{\kappa_{\text{SC}}^2(\mathbf{X})}^\alpha(\mu n^3v) \nonumber \\
    &= \frac{e^{-\frac{1}{\mu v}-\rho}}{\mu^{\alpha+1} v^{\alpha+2}} \sum_{k_1=0}^{\infty}\text{...}\sum_{k_\alpha=0}^{\infty} \left(\prod_{j=1}^\alpha \frac{1}{(j+k_j+1)!k_j!(\mu v)^{k_j}} \right) \lim_{n\xrightarrow{}\infty} \Theta(n,v)  H(v)\end{align}
    where 
    \begin{align}
    \Theta(n,v)= \det\left[\frac{(n+\alpha)!\mathcal{G}_i(\mu n^3v,\frac{\rho}{n})}{n^{3\alpha+2\sum_{j=1}^\alpha k_j+3}}\;\; \prod_{\ell=0}^{\alpha-i-1}(\Tilde{c_j}-\ell)\right]_{\substack{i=0,...,\alpha \\ j=1,...,\alpha}} .
 \end{align}
 To facilitate further analysis, we make use of the elementary row operations shown in \cite[Lemma A.1]{ref:prathapSIAM} on $\Theta(n,v)$ and expand the resultant determinant using its first column to obtain
  \begin{align}
 \label{eq_asymp_det1}
     \Theta(n,v) &= \Delta_\alpha(\textbf{c})\sum_{i=0}^\alpha \sum_{\ell=0}^{\alpha-i} \sum_{k=0}^i (-1)^i S_{\alpha-i}^{(\alpha-i-\ell)}  \frac{(n+\alpha)!\mathcal{G}_{i+\ell}(\mu n^3v,\frac{\rho}{n})}{n^{3\alpha+2\sum_{j=1}^\alpha k_j+3}} a_{k,i} n^k
 \end{align}{}
 where $\textbf{c}=\left\{c_1,...,c_\alpha\right\}$ with $c_j=j+k_j$, $a_{k,i}$'s are constant coefficients independent of $n$ with $a_{0,0}=1$, and $S_n^{(m)}$ is the Stirling number of the second kind with $S_\alpha^{(\alpha)}=1$ \cite{ref:gradshteyn}. Noting the fact that
 \begin{align}
     \lim_{n\xrightarrow{}\infty}(-1)^i S_{\alpha-i}^{(\alpha-i-\ell)} \frac{(n+\alpha)!\mathcal{G}_{i+\ell}(\mu n^3v, \frac{\rho}{n})}{n^{3\alpha+2\sum_{j=1}^\alpha k_j+3}} a_{k,i} n^k &= \left\{ \begin{matrix}{} e^\rho & \text{ for } i,\ell=0 \\ 0 & \text{ otherwise } \end{matrix} \right.,
 \end{align}
we conclude $\displaystyle\lim_{n\xrightarrow{}\infty} \Theta(n,v)=e^\rho \Delta_\alpha(\textbf{c})$, which upon substituting into (\ref{eq_asymp_substitution1}) gives
 \begin{align}
    f_V^\alpha(v) &= \frac{e^{-\frac{1}{\mu v}}}{\mu^{\alpha+1} v^{\alpha+2}} \sum_{k_1=0}^{\infty}\text{...}\sum_{k_\alpha=0}^{\infty} \left(\prod_{j=1}^\alpha \frac{1}{(j+k_j+1)!k_j! (\mu v)^{k_j}} \right) \Delta_\alpha(\textbf{c}) H(v).
 \end{align}
 Finally, following the developments in \cite{ref:prathapJMVA}, we obtain (\ref{asy pdf theo}) which concludes the proof. }
\end{IEEEproof}
The above limiting p.d.f. is significantly less complicated than that of the analytical expression which can be obtained by taking the derivative of the limiting c.d.f. given in (\ref{eq dem asym}). It is noteworthy that exactly the same limiting c.d.f. has been obtained in \cite{ref:prathapJMVA} when $\mathbf{X}$ is distributed as uncorrelated complex Gaussian with rank-one mean and in \cite{ref:akemann} for $\mathbf{X}$ having independent complex standard normal entries. 
Since the above limiting p.d.f. is independent of $\eta$, we can expect it to coincide with the limiting p.d.f. corresponding to $\eta=0$ case  given in \cite[Eq. 4.1]{ref:prathapSIAM}. Although the two expressions look different, a sanity check reveals that they are two different representations of the same analytical expression.

\begin{figure}[tbhp]
    \centering
    \includegraphics[width=12cm]{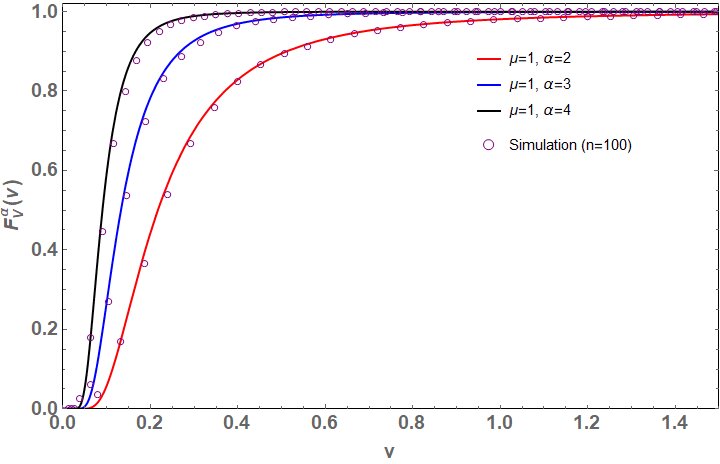}
    \caption{Comparison of simulated data points and the analytical asymptotic c.d.f. $F^\alpha_V(v)$ for different $\alpha$.}
    \label{fig_asymp}
\end{figure}

\begin{figure}[tbhp]
    \centering
   \includegraphics[width=12cm]{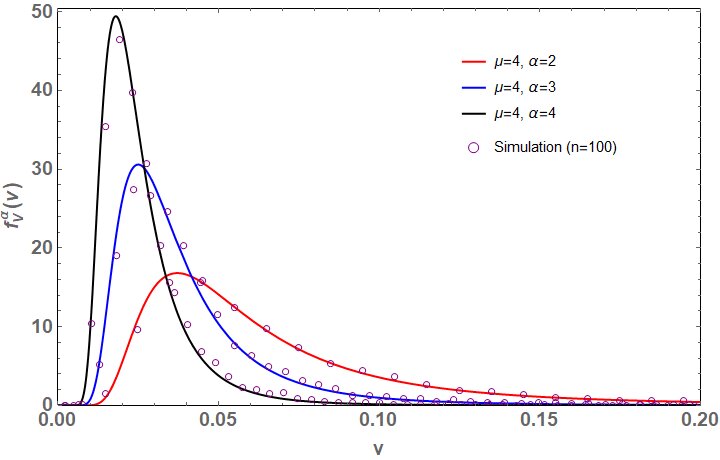}
    \caption{Comparison of simulated data points and the analytical asymptotic p.d.f. $f^\alpha_V(v)$ for different $\alpha$.}
    \label{fig_asymp_pdf}
\end{figure}

Figure \ref{fig_asymp} compares the analytical asymptotic c.d.f. given by Corollary \ref{cor asy cdf} with the simulated data points. This figure further highlights the advantage of the asymptotic formula since it compares favourably with finite $n$ results. Moreover, Fig. \ref{fig_asymp_pdf} depicts the comparison between the  theoretical limiting p.d.f. expression given by (\ref{asy pdf theo}) and corresponding simulation results. Both figures reveal the accuracy of our limiting expressions, particularly at the tail of the p.d.f. even for finite values of $n$. 

\begin{figure}[tbhp]
    \centering
   \includegraphics[width=12cm]{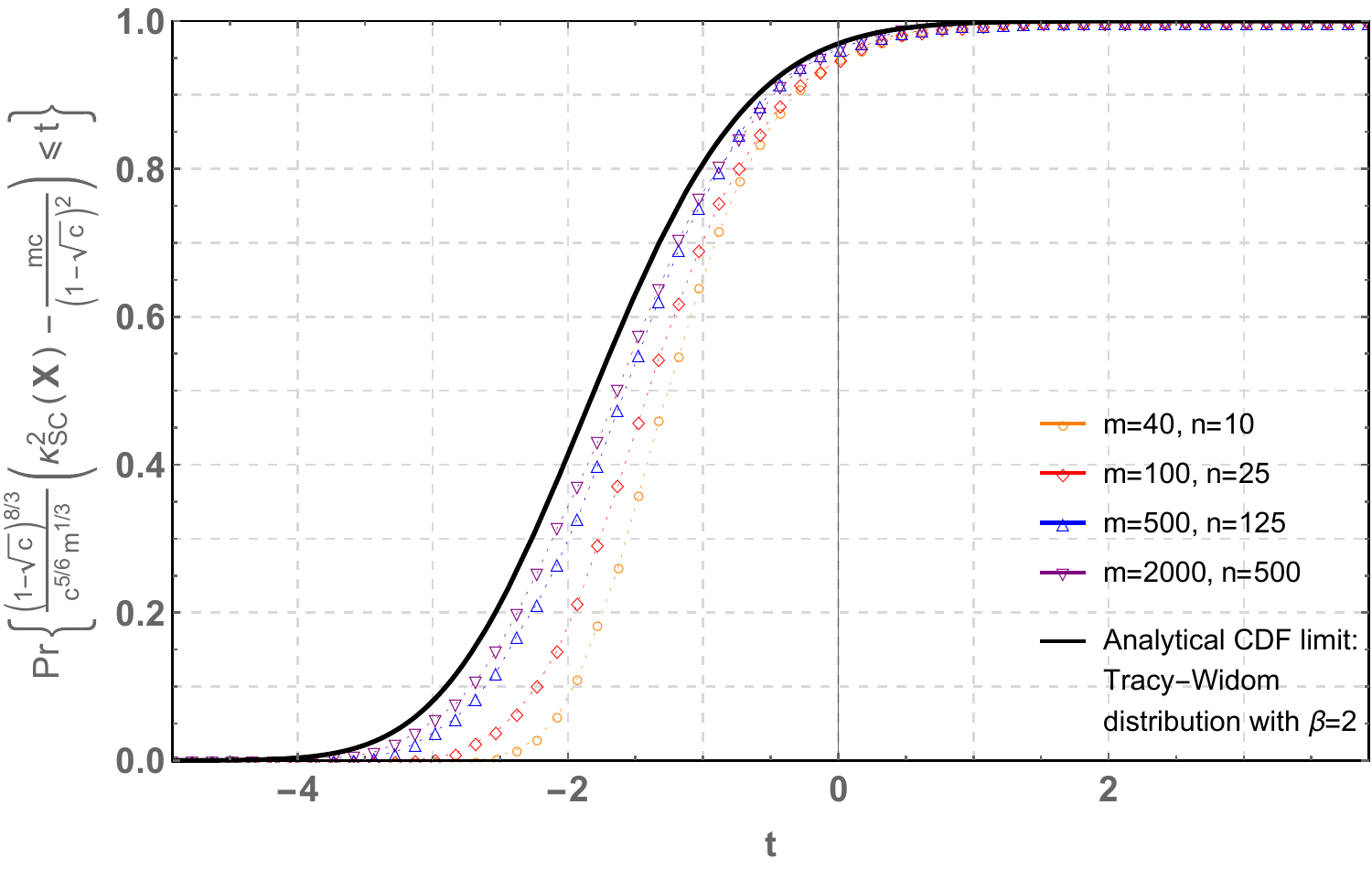}
    \caption{ The c.d.f. of $m^{-\frac{1}{3}}\left(1-\sqrt{c}\right)^{\frac{8}{3}}
    c^{-\frac{5}{6}}\left(\kappa_{\text{SC}}^{2}(\mathbf{X})-\frac{mc}{(1-\sqrt{c})^2}\right)$ for different $m,n$ configurations with $\eta=1$ and $c=0.25$. The limiting Tracy-Widom distribution corresponding to $\beta=2$ (i.e., complex case) is also shown for comparison.}
    \label{fig_tw_law_1}
\end{figure}

\begin{figure}[tbhp]
    \centering
   \includegraphics[width=12cm]{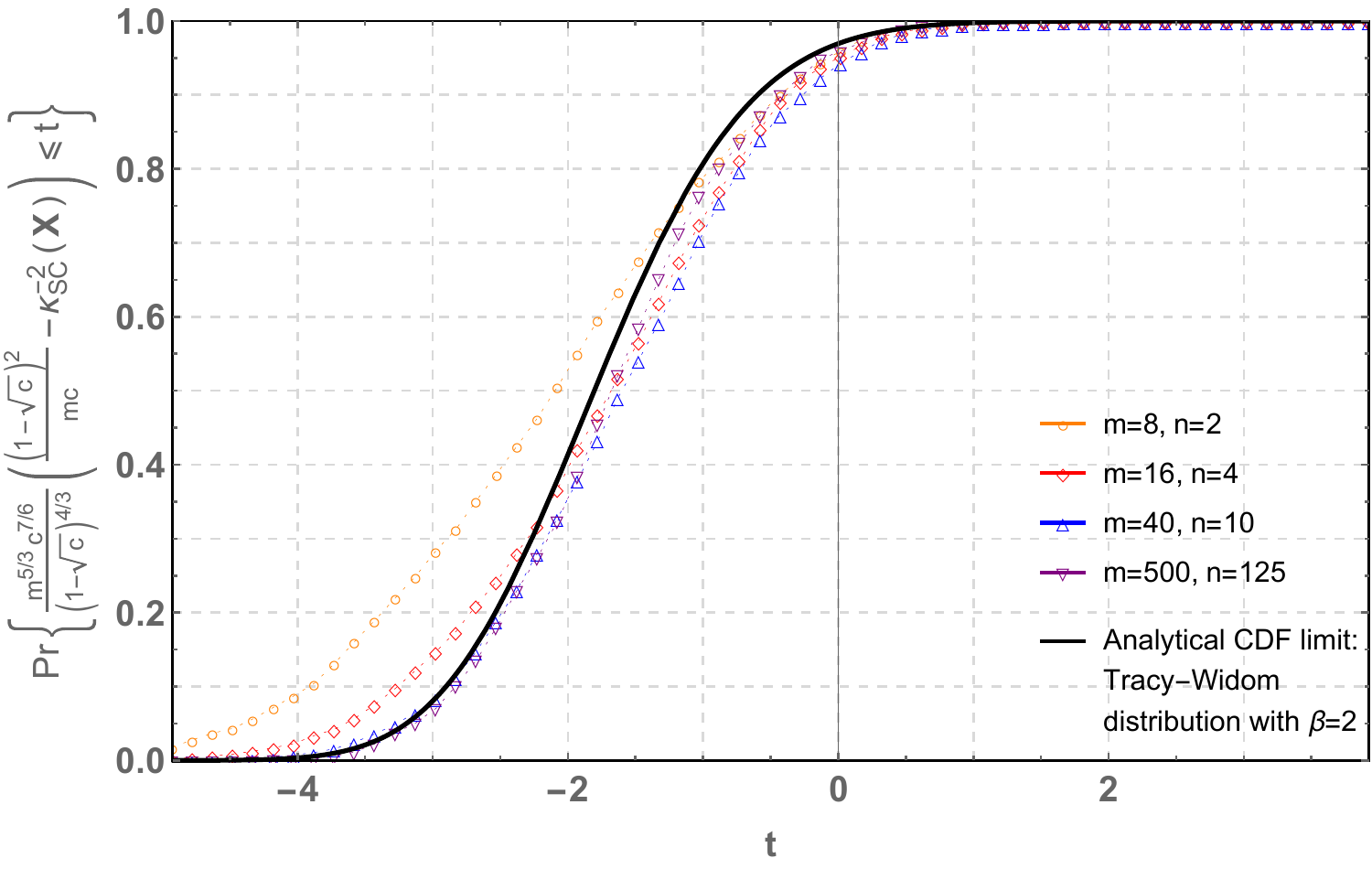}
    \caption{The c.d.f. of $m^{\frac{5}{3}}\left(1-\sqrt{c}\right)^{-\frac{4}{3}}c^{\frac{7}{6}}
    \left(\frac{(1-\sqrt{c})^2}{mc}-\kappa_{\text{SC}}^{-2}(\mathbf{X})\right)$ for different $m,n$ configurations with $\eta=1$ and $c=0.25$. The limiting Tracy-Widom distribution corresponding to $\beta=2$ (i.e., complex case) is also shown for comparison.}
    \label{fig_tw_law_2}
\end{figure}

{\color{blue}
\subsection{Limiting Distribution in the $n/m\to c\in (0,1)$ Regime}

Here we derive the limiting distribution of suitably centered and scaled $\kappa^2_{\rm SC}(\mathbf{X})$ as $m,n\to\infty$ such that $n/m\to c\in(0,1)$. 

Since $\kappa^2_{\text{SC}}(\mathbf{X})$ is intimately related to the behavior of $\lambda_1$ (i.e., the minimum eigenvalue) in this regime, we need to focus on the limiting distribution of $\lambda_1$.
In this respect,
as $m,n\to\infty$ such that $n/m\to c\in(0,1)$, for correlated Wishart matrices with their covariance matrices having a limiting spectral measure\footnote{Let the eigenvalues of the covariance matrix $\boldsymbol{\Sigma}$ be $0<\varrho_1\leq \varrho_2\leq\ldots\leq \varrho_n$. Then we assume that the spectral measure $\nu_n=\frac{1}{n}\sum_{k=1}^n \delta_{\varrho_k}$, where $\delta_x$ denotes the Direc measure at point $x$, converges weakly towards a limiting distribution $\nu$ as $m,n\to\infty$ such that $n/m\to c\in(0,1)$ \cite{ref:walid}.},
a Tracy-Widom type convergence result has been established in \cite{ref:walid}. To be specific, following \cite[Theorem 3]{ref:walid} and noting that $\nu({\rm d}x)=\delta_1 {\rm d}x$ for the single-spiked model (i.e., $\boldsymbol{\Sigma}=\mathbf{I}_n+\eta \mathbf{uu}^*$), we may write the limiting distribution of $\lambda_1$, after some tedious algebraic manipulation, as
\begin{align}
\label{mindoubleasy}
   \lim_{m\to\infty} \Pr\left\{\frac{1}{m^{\frac{1}{3}}}
    \frac{c^{\frac{1}{6}}}{\left(1-\sqrt{c}\right)^{\frac{4}{3}}}\left(m(1-\sqrt{c})^2-\lambda_1\right)\leq t\right\}=F_2(t)
\end{align}
where $F_2(t)$ denotes the famous Tracy-Widom distribution \cite{ref:tracy} corresponding to $\beta=2$ (i.e., complex case)
\begin{align}
    F_2(t)=\exp\left(-\int_t^\infty (x-t) q^2(x) {\rm d}x\right)
\end{align}
in which $q(x)$ denotes the Hastings-McLeod solution of the homogeneous Painlev\'{e} II equation $\frac{{\rm d^2}}{{\rm d} x^2} q(x)=2 q^3(x)+xq(x)$ characterized by the boundary condition $q(x)\sim {\rm Ai}(x)$ as $x\to\infty$ with ${\rm Ai}(x)$ denoting the Airy function\footnote{The Airy function is characterized in turn by $\frac{{\rm d^2}}{{\rm d} x^2} {\rm Ai}(x)=x{\rm Ai}(x)$ and $ {\rm Ai}(+\infty)=0$ \cite{ref:forresterLogGases}.}.
Since we are interested in the asymptotic characterization of $\kappa_{\text{SC}}^{2}(\mathbf{X})$, we consider the centered and scaled random variable $\frac{\left(1-\sqrt{c}\right)^{\frac{8}{3}}}{c^{\frac{5}{6}}m^{\frac{1}{3}}}\left(\kappa_{\text{SC}}^{2}(\mathbf{X})-\frac{mc}{(1-\sqrt{c})^2}\right)$, which can be rewritten as
\begin{align*}
 \frac{\left(1-\sqrt{c}\right)^{\frac{8}{3}}}{c^{\frac{5}{6}}m^{\frac{1}{3}}}\left(\kappa_{\text{SC}}^{2}(\mathbf{X})-\frac{mc}{(1-\sqrt{c})^2}\right)
 &=
 \frac{m^{\frac{2}{3}}\left(1-\sqrt{c}\right)^{\frac{8}{3}} }{c^{\frac{5}{6}}}
 \left(\frac{\sum_{k=1}^n \lambda_k/m^2}{\lambda_1/m}-\frac{c}{(1-\sqrt{c})^2}\right)\nonumber\\
 &= m^{\frac{2}{3}}\left(1-\sqrt{c}\right)^{\frac{2}{3}} c^{\frac{1}{6}}
 \frac{\left(c^{-1}(1-\sqrt{c})^2 \sum_{k=1}^n \lambda_k/m^2- \lambda_1/m\right)}{ \lambda_1/m}.
\end{align*}
Consequently,  noting that $\lambda_1/m$ converges almost surely to $(1-\sqrt{c})^2$ \cite{ref:baikPhaseTrans} and $\sum_{k=1}^n \lambda_k/m^2$ converges in probability to $c$\footnote{It can be proved that the characteristic function of $\sum_{j=1}^n\lambda_j/m^2=\text{tr}(\mathbf{W})/m^2$ given by $\mathbb{E}\left\{e^{-j\omega\text{tr}(\mathbf{W})/m^2 }\right\}=\frac{1}{\left(1+(1+\eta)j\omega/m^2\right)^{m} (1+j\omega/m^2)^{m(n-1)}}$ converges to $e^{-j\omega c}$ as $m,n\to\infty$ such that $n/m\to c\in(0,1)$ when $\eta=O(1)$. Therefore, following the L\'evy's continuity theorem, we obtain $\text{tr}(\mathbf{W})/m^2$ converges to $c$ weakly (i.e., converges in distribution). The final claim follows by noting that weak convergence to a constant amounts to convergence in probability to the same constant.} given $\eta=O(1)$, we may use Slutsky's lemma \cite{ref:vaart}, in view of (\ref{mindoubleasy}), to obtain
\begin{align}
\label{scndoublimit}
  \lim_{m\to\infty}  \Pr\left\{
    \frac{\left(1-\sqrt{c}\right)^{\frac{8}{3}}}{c^{\frac{5}{6}}m^{\frac{1}{3}}}\left(\kappa_{\text{SC}}^{2}(\mathbf{X})-\frac{mc}{(1-\sqrt{c})^2}\right)\leq t\right\}&=F_2(t).
\end{align}
A similar procedure can be used to establish
\begin{align}
    \label{eq_tw_1}
   \lim_{m\to\infty}  \Pr\left\{
    \frac{m^{\frac{5}{3}}c^{\frac{7}{6}}}{\left(1-\sqrt{c}\right)^{\frac{4}{3}}}\left(\frac{(1-\sqrt{c})^2}{mc}-\kappa_{\text{SC}}^{-2}(\mathbf{X})\right)\leq t\right\}&=F_2(t).
\end{align}
A careful inspection of the above c.d.f.s revels that properly centred $\kappa_{\text{SC}}^{2}(\mathbf{X})$ fluctuates on the scale $m^{\frac{1}{3}}$, whereas $\kappa_{\text{SC}}^{-2}(\mathbf{X})$ fluctuates on the scale $m^{-\frac{5}{3}}$. Moreover, in the light of observation that, for $\boldsymbol{\Sigma}=\mathbf{I}_n$, we have  \cite{ref:walid}
\begin{align}
  \lim_{m\to\infty}  \Pr\left\{
    \frac{\left(1-\sqrt{c}\right)^{\frac{8}{3}}}{c^{\frac{5}{6}}m^{\frac{1}{3}}}\left(\kappa_{\text{SC}}^{2}(\mathbf{X})-\frac{mc}{(1-\sqrt{c})^2}\right)\leq t\right\}&=F_2(t),
\end{align}
which is identically equal to (\ref{scndoublimit}),
we conclude that $\kappa_{\text{SC}}^{2}(\mathbf{X})$ does not have discrimination power to detect the presence of a weak signal. In other words, related to the CR spectrum sensing, when $\eta_{\rm snr}=\gamma||\mathbf{h}||^2/\sigma^2=O(1)$ (i.e., low SNR regime), the probability of detection of the test statistic $T$ converges to zero as $m,n\to\infty$ such that $n/m\to c\in(0,1)$. However, a sanity check reveals that, for $\eta_{\rm snr}=O(n)$ (i.e., high SNR regime), the test statistic $T$ still retains its detection power in the above asymptotic regime.

To further highlight the above asymptotic behavior, in Figs. \ref{fig_tw_law_1} and \ref{fig_tw_law_2}, we plot the c.d.f.s of $m^{-\frac{1}{3}}\left(1-\sqrt{c}\right)^{\frac{8}{3}}
    c^{-\frac{5}{6}}\left(\kappa_{\text{SC}}^{2}(\mathbf{X})-\frac{mc}{(1-\sqrt{c})^2}\right)$ and $m^{\frac{5}{3}}\left(1-\sqrt{c}\right)^{-\frac{4}{3}}c^{\frac{7}{6}}
    \left(\frac{(1-\sqrt{c})^2}{mc}-\kappa_{\text{SC}}^{-2}(\mathbf{X})\right)$, respectively. The results are shown for different $m,n$ configurations with $c=0.25$ and $\eta=1$. The limiting Tracy-Widom distribution corresponding to $\beta=2$ (i.e., complex case) is also shown for comparison. As can be seen from the figures, the limiting c.d.f.s do not compare favourably with finite dimensional results; particularly, for small $m,n$ configurations. Nevertheless, as $m$ and $n$ diverge, the Tracy-Widom c.d.f. serves as a good approximation.}

\section{Conclusions}
This paper investigates the exact p.d.f. characterization of $\kappa_{\text{SC}}^2(\mathbf{X})$ for row correlated complex Gaussian distributed $\mathbf{X}$ with the so called single-spiked covariance matrix. Capitalizing on the powerful orthogonal polynomial approach from finite dimensional random matrix theory, we derive an exact expression for the p.d.f. which contains the determinant of a square matrix whose dimensions depend on the relative difference $m-n$ (i.e., the degree of {\it rectangularity} of $\mathbf{X}$). {\color{blue}To demonstrate the utility of this new expression, noting the significance of $\kappa_{\text{SC}}^2(\mathbf{X})$ as one of the test statistics used in the CR blind spectrum sensing, we develop corresponding ROC curves in various scenarios. It turns out that, when $m=n$, an explicit analytical relationship between the probabilities of detection and false alarm can be obtained. Nevertheless, obtaining such an analytical relationship for $m\neq n$ scenario seems intractable.}
Subsequently, we characterize the behavior of the scaled $\kappa_{\text{SC}}^2(\mathbf{X})$ in the asymptotic regime, where $m,n\to\infty$ with $m-n$ fixed, by deriving the limiting c.d.f. and p.d.f. expressions. In particular, our analytical stochastic convergence results reveal that, if $\eta$ scales on the order of $1/n$, then $\kappa_{\text{SC}}^2(\mathbf{X})$ scales on the order of $n^3$ in the former asymptotic regime. Although derived for asymptotically large $m$ and $n$, these limiting distributions compare favourably with finite $m$ and $n$ results as well. {\color{blue}On the other hand, we also establish a Tracy-Widom class of stochastic convergence result for  $\kappa_{\text{SC}}^2(\mathbf{X})$ as  $m,n\to\infty$ such that $n/m\to c\in(0,1)$. In this respect, we show that properly centered $\kappa_{\text{SC}}^2(\mathbf{X})$ fluctuates on the scale $m^{1/3}$. Nevertheless, our analytical results reveal that $\kappa_{\text{SC}}^2(\mathbf{X})$ (also $T(\boldsymbol{\lambda})$) does not have discrimination power to detect a weak signal (i.e., when $\eta=O(1)$ or equivalently $\eta_{\rm snr}=O(1)$) in this asymptotic regime.}  

The next natural question is whether we can employ the same analytical machinery to extend the above results to the case in which the covariance matrix has an arbitrary number of spikes. An affirmative answer in this respect requires further research and therefore, remains as an open problem.

\appendices
\section{Proof of Corollary \ref{corretazero}}
\label{appendix_corretazero}
By substituting $\eta=0$ in (\ref{eq_pdf_Z}), we get
\begin{multline}
    f_{\kappa_{\text{SC}}^2(\mathbf{X})}^\alpha(z) = \frac{(n+\alpha)!H(z-n)}{z^{n(n+\alpha)}} \sum_{k_1=0}^{n+\alpha-2} \text{...} \sum_{k_{\alpha}=0}^{n-1}\left(\prod_{j=1}^{\alpha}\frac{(n+\alpha-j-1)!}{(j+k_j+1)!k_j!} \right) \\ 
    \times
    (z-n)^{(n-1)(n+\alpha+1)-\sum_{j=1}^{\alpha}k_j-1} \det\left[ G_i(z) \;\; \frac{1}{\Gamma(n+i-j-k_j)} \right]_{\substack{i=0,..,\alpha \\ j=1,..,\alpha}}
\end{multline}
where
\begin{align}
    G_i(z) &= \left\{ \begin{matrix}G_0(z) &\text{ for }i = 0 \\ 0 &\text{ otherwise } \end{matrix} \right.,
\end{align}
with
\begin{align*}
    G_0(z)
    &= \frac{\Gamma\left(n(n+\alpha)\right)}{\Gamma(n)\Gamma\left((n-1)(n+\alpha+1)-\sum_{j=1}^{\alpha}k_j\right)} \text{ .}
\end{align*}
Now we expand the determinant using its first column to 
yield
\begin{multline}
    f_{\kappa_{\text{SC}}^2(\mathbf{X})}^\alpha(z) = \frac{(n+\alpha)!H(z-n)}{z^{n(n+\alpha)}} \sum_{k_1=0}^{n+\alpha-2} \text{...} \sum_{k_{\alpha}=0}^{n-1}\left(\prod_{j=1}^{\alpha}\frac{(n+\alpha-j-1)!}{(j+k_j+1)!k_j!} \right) \\ 
    \times
    (z-n)^{(n-1)(n+\alpha+1)-\sum_{j=1}^{\alpha}k_j-1} \frac{\Gamma\left(n(n+\alpha)\right)}{\Gamma(n)\Gamma\left((n-1)(n+\alpha+1)-\sum_{j=1}^{\alpha}k_j\right)} \\ 
    \times \det\left[\frac{1}{\Gamma(n+i-j-k_j)} \right]_{\substack{i,j=1,..,\alpha}}.
\end{multline}
To facilitate further analysis, let us rearrange the terms in the determinant to obtain
\begin{multline}
    f_{\kappa_{\text{SC}}^2(\mathbf{X})}^\alpha(z) = 
    \Gamma\left(n(n+\alpha)\right)\left(\prod_{j=0}^\alpha n+j\right)(z-n)^{n(n+\alpha)-\alpha-2}z^{-n(n+\alpha)}\left(\prod_{j=1}^\alpha\frac{1}{(j+1)!}\right) \\
    \times \sum_{k_1=0}^{n+\alpha-2} \text{...} \sum_{k_{\alpha}=0}^{n-1}\left(\prod_{j=1}^{\alpha}\frac{(n+\alpha-j-1)!(j+1)!}{(j+k_j+1)!k_j!(n+\alpha-j-1-k_j)!} (z-n)^{-k_j}\right) \\
    \times \frac{H(z-n)}{\Gamma\left(n(n+\alpha)-\alpha-1-\sum_{j=1}^{\alpha}k_j\right)}\det\left[\frac{(n+\alpha-j-1-k_j)!}{(n+i-j-k_j-1)!} \right]_{\substack{i,j=1,..,\alpha}}
\end{multline}
which simplifies after some algebraic manipulations giving
\begin{multline}
    f_{\kappa_{\text{SC}}^2(\mathbf{X})}^\alpha(z) = \left(n(n+\alpha)-1\right)!\left(\prod_{j=0}^\alpha\frac{n+j}{(j+1)!}\right)(z-n)^{n(n+\alpha)-\alpha-2}z^{-n(n+\alpha)} \\
    \times \sum_{k_1=0}^{n+\alpha-2} \text{...} \sum_{k_{\alpha}=0}^{n-1}\left(\prod_{j=1}^{\alpha}(-1)^{k_j}\frac{(-n-\alpha+j+1)_{k_j}}{(j+1)_{k_j}k_j!} (z-n)^{-k_j}\right) \\
    \times \frac{H(z-n)}{\Gamma\left(n(n+\alpha)-\alpha-1-\sum_{j=1}^{\alpha}k_j\right)}\det\left[\prod_{\ell=0}^{\alpha-i-1}(n+\alpha-j-k_j-1-\ell) \right]_{\substack{i,j=1,..,\alpha}}.
\end{multline}
Finally, we invoke \cite[Lemma A.1]{ref:prathapSIAM} to further simplify the above determinant and apply some algebraic manipulations to obtain (\ref{eq dem pdf white}), which concludes the proof.

\section{Derivation of the p.d.f. of the minimum eigenvalue}
\label{appendix_pdfmineig}
Let the p.d.f be $f_{\lambda_1}^\alpha(x)$. By integrating the joint p.d.f of the eigenvalues given in (\ref{cor joint}) in terms of $\lambda_2,...,\lambda_n$, we get
\begin{align}
    f_{\lambda_1}^\alpha(x) &= \int_{\mathcal{R}} f(x,\lambda_2,\dotsc,\lambda_n) {\rm d}\lambda_2 \dotsc {\rm d}\lambda_n
\end{align}
where $\mathcal{R}=\{x\le \lambda_2 \le \dotsc \le \lambda_n \}$. Noting that the integral is symmetric with respect to $\lambda_2,\dotsc,\lambda_n$, we change the ordered region of integration into an $(n-1)$-fold set of unordered regions, to obtain
\begin{multline}
    \displaystyle
    f_{\lambda_1}^\alpha(x) = \frac{C_{n, \alpha, \eta}}{(n-1)!} \int_{\left[x, \infty\right)^{n-1}} x^\alpha e^{-x} \prod_{i=2}^n\lambda_i^\alpha e^{-\lambda_i} \prod_{i=2}^n\left(\lambda_i-x\right)^2 \Delta_{n-1}^2(\boldsymbol{\lambda})  \\
    \times \left( \frac{e^{c_\eta x}}{\displaystyle\prod_{i=2}^n\left (x-\lambda_i\right)} + \sum_{k=2}^n\frac{e^{c_\eta \lambda_k}}{\left(\lambda_k-x\right)\displaystyle\prod_{\substack{i=2 \\ i\neq k}}^n\left(\lambda_k-\lambda_i\right)} \right) {\rm d}\lambda_2 \dotsc {\rm d}\lambda_n \text{.}
\end{multline}
For further simplification, we split the above integral  to yield
\begin{align}
\label{min pdf ab decom}
    f_{\lambda_1}^\alpha(x) &= \tilde{\mathcal{A}}(x) + \tilde{\mathcal{B}}(x)
\end{align}
where
\begin{equation}
\label{pdf min a def}
    \tilde{\mathcal{A}}(x) = \frac{C_{n, \alpha, \eta}}{(n-1)!} x^\alpha e^{-x(1-c_\eta)} \int_{\left[x, \infty\right)^{n-1}} \prod_{i=2}^n\lambda_i^\alpha e^{-\lambda_i} \left(x-\lambda_i\right)  \Delta_{n-1}^2(\boldsymbol{\lambda}) {\rm d}\lambda_2 \dotsc {\rm d}\lambda_n
\end{equation}
and
\begin{multline}
\label{pdf min b def}
    \tilde{\mathcal{B}}(x) = \frac{C_{n, \alpha, \eta}}{(n-1)!} x^\alpha e^{-x} \int_{\left[x, \infty\right)^{n-1}} \sum_{k=2}^n\frac{e^{c_\eta \lambda_k}}{\left(\lambda_k-x\right)\displaystyle\prod_{\substack{i=2 \\ i\neq k}} ^ n\left(\lambda_k-\lambda_i\right)}\\
    \times \prod_{i=2}^n\lambda_i^\alpha e^{-\lambda_i} \left(\lambda_i-x\right)^2 \Delta_{n-1}^2(\boldsymbol{\lambda})  {\rm d}\lambda_2 \dotsc {\rm d}\lambda_n \text{.}
\end{multline}
Let us now focus on simplifying $\tilde{\mathcal{A}}(x)$. To this end, we apply the variable transformations, $\lambda_i-x=y_{i-1}$,  $i=2,\dotsc,n$, to (\ref{pdf min a def}) with some algebraic manipulations to obtain 
\begin{align}
    \tilde{\mathcal{A}}(x) &=
    \frac{(-1)^{(\alpha+1)(n-1)}}{(n-1)!}C_{n, \alpha, \eta}x^\alpha e^{-x(n-c_\eta)} R_{n-1}^{(\alpha)}(-x)
\end{align}
from which we obtain in view of (\ref{eq_R_eval})
\begin{multline}
\label{eq a sol}
    \tilde{\mathcal{A}}(x) = \frac{(-1)^{n-1}}{(n-1)!} C_{n, \alpha, \eta}x^\alpha e^{-x(n-c_\eta)} \prod_{j=0}^{n-2}(j+1)!(j+1)! \prod_{j=0}^{\alpha-1}\frac{(n+j-1)!}{j!}\\
    \times \det \left[L_{n+i-j-1}^{(2)}\right]_{\substack{i,j=1\dotsc \alpha}} \text{.}
\end{multline}
To further simplify $\tilde{\mathcal{B}}(x)$, noting that due to symmetry, each term in the sum contributes the same amount to the total, we rewrite (\ref{pdf min b def}) after some algebraic manipulations as
\begin{multline}
    \tilde{\mathcal{B}}(x) =  \frac{C_{n, \alpha, \eta}x^\alpha e^{-x}}{(n-2)!} \int_x^\infty \int_{\left[x,\infty\right)^{n-2}} \lambda_2^\alpha e^{-\lambda_2(1-c_\eta)}(\lambda_2-x) \prod_{i=3}^n \lambda_i^\alpha e^{-\lambda_i}(\lambda_i-x)^2(\lambda_2-\lambda_i)\\
    \times \Delta_{n-2}^2(\boldsymbol{\lambda})  {\rm d}\lambda_3 \dotsc {\rm d}\lambda_n {\rm d}\lambda_2.
\end{multline}
Now it is convenient to introduce the variable transformations, $\lambda_2-x=y$ and $\lambda_i-x=y_{i-2}$, $i=3,\dotsc,n$, to yield
\begin{align}
    \tilde{\mathcal{B}}(x) 
    &= \frac{(-1)^{n\alpha}}{(n-2)!}C_{n,\alpha,\eta}x^\alpha e^{-x(n-c_\eta)} \int_0^\infty (y+x)^\alpha e^{-y(1-c_\eta)}y T_{n-2}^{(\alpha)}(y, -x) {\rm d}y,
\end{align}
from which we obtain in view of  (\ref{Teval})
\begin{multline}
    \tilde{\mathcal{B}}(x)= \frac{(-1)^{n}}{(n-2)!}C_{n,\alpha,\eta}\overline{\mathcal{K}}_{n-2, \alpha} x^\alpha e^{-x(n-c_\eta)} \int_0^\infty e^{-y(1-c_\eta)} y  \\
    \times \det\left[L_{n+i-3}^{(2)}(y) \;\; L_{n+i-j-1}^{(j)}(-x)\right]_{\substack{i=1,\dotsc,\alpha+1 \\ j=2,\dotsc,\alpha+1}} {\rm d}y \text{.}
\end{multline}
Since only the first column of the above determinant depends on $y$, we rewrite the above integral as
\begin{align}
\label{eq_Bx_with_taui}
    \tilde{\mathcal{B}}(x) &= \frac{(-1)^{n}}{(n-2)!}C_{n,\alpha,\eta}\overline{\mathcal{K}}_{n-2, \alpha} x^\alpha e^{-x(n-c_\eta)} \det\left[\tau_i \;\; L_{n+i-j-1}^{(j)}(-x)\right]_{\substack{i=1,\dotsc,\alpha+1 \\ j=2,\dotsc,\alpha+1}}
\end{align}
where
\begin{align}
\label{eq_taui_integral}
    \tau_i &= \int_0^\infty y e^{-y(1-c_\eta)} L_{n+i-3}^{(2)}(y) {\rm d}y \text{.}
\end{align}
We may use the contiguity relation (\ref{lagcont}) and \cite[Eq. 7.414.5]{ref:gradshteyn} with some tedious algebraic manipulations to obtain
\begin{align}
    \tau_i &=1+\gamma_i
\end{align}
where $\gamma_i=(-1)^{n}(\eta+1)\eta^{n-1}(n+i-1-c_\eta)(-\eta)^{i-1}$.
Substituting this result back into (\ref{eq_Bx_with_taui}) and expanding the resultant determinant using the multi-linear property gives us
\begin{multline}
\label{eq_Bx_two_dets}
    \tilde{\mathcal{B}}(x) 
    = \frac{(-1)^{n}}{(n-2)!}C_{n,\alpha,\eta}\overline{\mathcal{K}}_{n-2, \alpha} x^\alpha e^{-x(n-c_\eta)} \det\left[1 \;\; L_{n+i-j-1}^{(j)}(-x)\right]_{\substack{i=1,\dotsc,\alpha+1 \\ j=2,\dotsc,\alpha+1}}  \\
    +\frac{(-1)^{n}}{(n-2)!}C_{n,\alpha,\eta}\overline{\mathcal{K}}_{n-2, \alpha} x^\alpha e^{-x(n-c_\eta)} \det\left[\gamma_i \;\; L_{n+i-j-1}^{(j)}(-x)\right]_{\substack{i=1,\dotsc,\alpha+1 \\ j=2,\dotsc,\alpha+1}}\text{.}
\end{multline}
Noting that (see (\ref{eq Lag det simp})) 
\begin{align*}
    \det\left[1 \;\; L_{n+i-j-1}^{(j)}(-x)\right]_{\substack{i=1,\dotsc,\alpha+1 \\ j=2,\dotsc,\alpha+1}} &= \det\left[L_{n+i-j-1}^{(j)}(-x)\right]_{\substack{i,j=1,\dotsc,\alpha}}
\end{align*}
and in view of (\ref{eq a sol}) we rewrite (\ref{eq_Bx_two_dets}) as
\begin{multline}
    \tilde{\mathcal{B}}(x) 
    \label{min pdf b last}
    = -\tilde{\mathcal{A}}(x)+\frac{(-1)^{n}
    }{(n-2)!}C_{n,\alpha,\eta}\overline{\mathcal{K}}_{n-2, \alpha} x^\alpha e^{-x(n-c_\eta)} \det\left[\gamma_i \;\; L_{n+i-j-1}^{(j)}(-x)\right]_{\substack{i=1,\dotsc,\alpha+1 \\ j=2,\dotsc,\alpha+1}}.
\end{multline}
Consequently, we use (\ref{min pdf b last}) in (\ref{min pdf ab decom}) with some algebraic manipulations to arrive at 
\begin{multline}
    f_{\lambda_1}^\alpha(x) 
    = \frac{(n-1)!x^\alpha e^{-x(n-c_\eta)}}{(n+\alpha-1)!(\eta+1)^\alpha} 
    \det\left[(n+i-1-c_\eta)(-\eta)^{i-1} \;\; L_{n+i-j-1}^{(j)}(-x)\right]_{\substack{i=1,\dotsc,\alpha+1 \\ j=2,\dotsc,\alpha+1}} \text{ .}
\end{multline}
Finally, 
we make use of the index shift, $i\rightarrow i-1, j\rightarrow j-1$, to conclude the proof.


\section{Derivation of the c.d.f. of the minimum eigenvalue}
\label{appendix_cdfmineig}
By definition, the c.d.f. of $\lambda_1$ can be written as
\begin{align}
    \label{eq_mineig_cdf_initial}
    F_{\lambda_1}^\alpha(x) &= \Pr(\lambda_1\leq x)= 1-\Pr(\lambda_1\ge x) \text{.}
\end{align}
As such, $\Pr(\lambda\ge x)$ can be written as 
\begin{align}
    \Pr(\lambda_1\ge x) &= \int_{x\le\lambda_1\le\dotsc\le\lambda_n}f(\lambda_1,\dotsc,\lambda_n) {\rm d}\lambda_1\dotsc {\rm d}\lambda_n \\
    &= \frac{C_{n, \alpha, \eta}}{n!} \int_{[x,\infty)^n} \prod_{i=1}^n \lambda_i^\alpha e^{-\lambda_i} \Delta_n^2(\boldsymbol{\lambda}) \sum_{k=1}^n \frac{e^{c_\eta \lambda_k}}{\displaystyle \prod_{\substack{i=1 \\ i\neq k}}^n (\lambda_k-\lambda_i)} {\rm d}\lambda_1\dotsc{\rm d}\lambda_n \label{eq cdf main der}
\end{align}
where we have exploited the fact that the integral is symmetric with respect to $\lambda_1\dotsc,\lambda_n$. Now it is worth observing that, due to symmetry, each term in the sum in the above expression contributes an equal amount to the total. Therefore, we can further simplify (\ref{eq cdf main der}) to yield   
\begin{equation}
    \Pr(\lambda_1\ge x) = \frac{C_{n,\alpha,\eta}}{(n-1)!}\int_{[x,\infty)^n} \lambda_1^\alpha e^{-\lambda_1(1-c_\eta)} \prod_{i=2}^n \lambda_i^\alpha e^{-\lambda_i}(\lambda_1-\lambda_i) \Delta_{n-1}^2(\boldsymbol{\lambda}) {\rm d}\lambda_1\dotsc {\rm d}\lambda_n \text{.}
\end{equation}
To facilitate further analysis, we introduce the variable transformations,
 $\lambda_1-x=y$ and $\lambda_i-x=y_{i-1}$ for $i=2,\dotsc,n$, to obtain
\begin{equation}
    \label{eq_prob_mineig_ge_x_with_Q}
    \Pr(\lambda_1\ge x)
    = \frac{(-1)^{(n-1)\alpha}}{(n-1)!}C_{n,\alpha,\eta}e^{-x(n-c_\eta)} \int_0^\infty (x+y)^\alpha e^{-y(1-c_\eta)} Q_{n-1}^{(\alpha)}(y, -x) {\rm d}y
\end{equation}
where
\begin{align*}
    Q_{n}^{(\alpha)}(a,b) &= \int_{[0,\infty)^n} \prod_{i=1}^n (a-y_i)(b-y_i)^\alpha e^{-y_i} \Delta_n^2(\boldsymbol{y}) {\rm d}y_1\dotsc{\rm d}y_n \text{.}
\end{align*}
Following \cite[Eq. A.4]{ref:prathapJMVA}, we obtain a closed-form solution to  $Q_n^{(\alpha)}(a,b)$ as
\begin{align}
    \label{eq_Q_eval}
    Q_n^{(\alpha)}(a,b) &= \frac{(-1)^{n+\alpha(n+\alpha)}{\hat{\mathcal{K}}}_{n, \alpha}}{(b-a)^\alpha} \det\left[L_{n+i-1}^{(0)}(a) \;\;\; L_{n+i+1-j}^{(j-2)}(b)\right]_{\substack{i=1,\dotsc,\alpha+1 \\ j=2,\dotsc,\alpha+1}}
\end{align}
where
\begin{align*}
    {\hat{\mathcal{K}}}_{n, \alpha} &= \frac{\prod_{j=1}^{\alpha+1}(n+j-1)!\prod_{j=0}^{n-1}j!(j+1)!}{\prod_{j=1}^{\alpha-1}j!}
\end{align*}
Consequently, we use (\ref{eq_Q_eval}) in (\ref{eq_prob_mineig_ge_x_with_Q}) to obtain
\begin{multline}
    \Pr(\lambda_1\ge x) = \frac{(-1)^{n\alpha}}{(n-1)!}C_{n,\alpha,\eta}e^{-x(n-c_\eta)}{\hat{\mathcal{K}}}_{n-1,\alpha} \int_0^\infty e^{-y(1-c_\eta)}  \\
    \times  \det\left[L_{n+i-2}^{(0)}(y) \;\;\; L_{n+i-j}^{(j-2)}(-x)\right]_{\substack{i=1,\dotsc,\alpha+1 \\ j=2,\dotsc,\alpha+1}} {\rm d}y \text{.}
\end{multline}
Observing that only the first column of the above determinant depends on $y$, we can further simplify the above integral to yield
\begin{equation}
\label{eq cdf min last}
    \Pr(\lambda_1\ge x) = \frac{(-1)^{n\alpha}}{(n-1)!}C_{n,\alpha,\eta}e^{-x(n-c_\eta)}{\hat{\mathcal{K}}}_{n-1,\alpha} \det\left[\nu_i \;\;\; L_{n+i-j}^{(j-2)}(-x)\right]_{\substack{i=1,\dotsc,\alpha+1 \\ j=2,\dotsc,\alpha+1}}
\end{equation}
where
\begin{equation}
    \nu_i = \int_0^\infty e^{-y(1-c_\eta)} L_{n+i-2}^{(0)}(y) {\rm d}y=(\eta+1)(-\eta)^{n+i-2}
\end{equation}
with the last equality follows from \cite[Eq. 7.414.6]{ref:gradshteyn}. Therefore, (\ref{eq cdf min last}) becomes, after some algebraic manipulations
\begin{equation}
    \Pr(\lambda_1\ge x)
    = \frac{e^{-x(1-c_\eta)}}{(\eta+1)^\alpha} \det\left[(-\eta)^{i-1} \;\;\; L_{n+i-j}^{(j-2)}(-x)\right]_{\substack{i=1,\dotsc,\alpha+1 \\ j=2,\dotsc,\alpha+1}} \text{.}
\end{equation}
Substituting this result into (\ref{eq_mineig_cdf_initial}) followed by the index shift, $i\rightarrow i-1, j\rightarrow j-1$, concludes the proof.


\ifCLASSOPTIONcaptionsoff
  \newpage
\fi



\bibliographystyle{IEEEtran}
\bibliography{./references.bib}
\end{document}